\documentclass[10pt]{article}
\usepackage{amsmath,amssymb}
\usepackage[dvips]{graphicx}
\usepackage{hyperref}
\newcommand{\status}{}

\newcommand{\file}{$\ti{\ }$/wisk/bran/bran.tex}
\renewcommand{\file}{}

\newcommand{\detail}[1]{\par\noi{\bf[Proof detail\ }{#1}
\hfill{\bf ]}\par\noi\hspace{-4pt}}
\renewcommand{\detail}[1]{}

\newcommand{\dis}{\displaystyle}
\newcommand{\txt}{\textstyle}

\newcommand{\noi}{\noindent}

\newcommand{\halmos}{\rule{1ex}{1.4ex}}
\def \qed {\nopagebreak{\hspace*{\fill}$\halmos$\medskip}}
\newcommand{\med}{\medskip}

\newtheorem{theorem}{Theorem}
\newtheorem{proposition}[theorem]{Proposition}
\newtheorem{corollary}[theorem]{Corollary}
\newtheorem{conjecture}[theorem]{Conjecture}
\newtheorem{lemma}[theorem]{Lemma}

\newtheorem{remark}[theorem]{Remark}

\newcommand{\bt}{\begin{theorem}}
\newcommand{\et}{\end{theorem}}
\newcommand{\bl}{\begin{lemma}}
\newcommand{\el}{\end{lemma}}
\newcommand{\bp}{\begin{proposition}}
\newcommand{\ep}{\end{proposition}}
\newcommand{\bcor}{\begin{corollary}}
\newcommand{\ecor}{\end{corollary}}
\newcommand{\br}{\begin{remark}\rm}
\newcommand{\er}{\end{remark}}
\newcommand{\bcon}{\begin{conjecture}}
\newcommand{\econ}{\end{conjecture}}

\newcommand{\be}{\begin{equation}}
\newcommand{\ee}{\end{equation}}
\newcommand{\bes}{\begin{equation*}}
\newcommand{\ees}{\end{equation*}}
\newcommand{\bee}{\begin{enumerate}}
\newcommand{\eee}{\end{enumerate}}
\newcommand{\bei}{\begin{itemize}}
\newcommand{\eei}{\end{itemize}}
\newcommand{\bea}{\begin{eqnarray}}
\newcommand{\eea}{\end{eqnarray}}
\newcommand{\beas}{\begin{eqnarray*}}
\newcommand{\eeas}{\end{eqnarray*}}
\newcommand{\ba}{\begin{array}}
\newcommand{\ea}{\end{array}}
\newcommand{\bc}{\be\begin{array}{r@{\,}c@{\,}l}}
\newcommand{\ec}{\end{array}\ee}

\newcommand{\al}{\alpha}
\newcommand{\bet}{\beta}
\newcommand{\ga}{\gamma}
\newcommand{\Ga}{\Gamma}
\newcommand{\de}{\delta}
\newcommand{\De}{\Delta}
\newcommand{\eps}{\varepsilon}
\newcommand{\la}{\lambda}
\newcommand{\La}{\Lambda}

\newcommand{\tet}{\theta}

\newcommand{\si}{\ensuremath{\sigma}}

\newcommand{\Ci}{{\cal C}}

\newcommand{\Ei}{{\cal E}}

\newcommand{\Gi}{{\cal G}}
\newcommand{\Hi}{{\cal H}}

\newcommand{\Li}{{\cal L}}
\newcommand{\Mi}{{\cal M}}
\newcommand{\Ni}{{\cal N}}

\newcommand{\Si}{{\cal S}}

\newcommand{\Ui}{{\cal U}}

\newcommand{\Xc}{{\cal X}}
\newcommand{\Yi}{{\cal Y}}

\newcommand{\R}{{\mathbb R}}
\newcommand{\N}{{\mathbb N}}
\newcommand{\Z}{{\mathbb Z}}

\newcommand{\E}{{\mathbb E}}
\renewcommand{\P}{{\mathbb P}}

\newcommand{\li}{\langle}
\newcommand{\re}{\rangle}

\newcommand{\up}{\uparrow}
\newcommand{\down}{\downarrow}
\newcommand{\sub}{\subset}
\newcommand{\beh}{\backslash}
\newcommand{\symdif}{\!\vartriangle\!}

\newcommand{\asto}[1]{\underset{{#1}\to\infty}{\longrightarrow}}
\newcommand{\Asto}[1]{\underset{{#1}\to\infty}{\Longrightarrow}}

\newcommand{\ti}{\tilde}
\newcommand{\dgg}{\dagger}
\newcommand{\ov}{\overline}
\newcommand{\un}{\underline}
\newcommand{\subb}[2]{_{\ba{c}\scriptstyle{#1}\\[-.15cm]\scriptstyle{#2}\ea}}

\newcommand{\ffrac}[2]{{\textstyle\frac{{#1}}{{#2}}}}
\newcommand{\dif}[1]{\ffrac{\partial}{\partial{#1}}}
\newcommand{\diff}[1]{\ffrac{\partial^2}{{\partial{#1}}^2}}
\newcommand{\difif}[2]{\ffrac{\partial^2}{\partial{#1}\partial{#2}}}

\newcommand{\di}{\mathrm{d}}
\newcommand{\half}{{[0,\infty)}}
\newcommand{\expo}{\mbox{\large\it e}}
\newcommand{\ex}[1]{\expo^{\,\textstyle{#1}}}

\newcommand{\Pois}{{\rm Pois}}
\newcommand{\Thin}{{\rm Thin}}
\newcommand{\quand}{\quad\mbox{and}\quad}
\newcommand{\var}{{\rm Var}}
\newcommand{\cov}{{\rm Cov}}
\makeatletter\@addtoreset{equation}{section}
\makeatother 

\begin{document}

\title{\vspace{-3cm}Systems of branching, annihilating, and coalescing
particles}
\author{
Siva R.~Athreya\vspace{6pt}\\
{\small Indian Statistical Institute}\\
{\small 8th Mile Mysore Road}\\
{\small RV College PO}\\
{\small Bangalore -560059, India\vspace{3pt}}\\
{\small e-mail: athreya@isibang.ac.in}\vspace{4pt}
\and
\and Jan M. Swart\vspace{6pt}\\
{\small \' UTIA}\\
{\small Pod vod\'arenskou v\v e\v z\' i 4}\\
{\small 18208 Praha 8}\\
{\small Czech Republic}\\
{\small e-mail: swart@utia.cas.cz}\vspace{4pt}}
\date{{\scriptsize\file}\quad\today}
\maketitle\vspace{-.7cm}
\status

\begin{abstract}\noi
This paper studies systems of particles following independent random walks and
subject to annihilation, binary branching, coalescence, and deaths. In the
case without annihilation, such systems have been studied in our 2005 paper
``Branching-coalescing particle systems''. The case with annihilation is
considerably more difficult, mainly as a consequence of the non-monotonicity
of such systems and a more complicated duality. Nevertheless, we
show that adding annihilation does not significantly change the long-time
behavior of the process and in fact, systems with annihilation can be obtained
by thinning systems without annihilation.
\end{abstract}

\vspace{.4cm}
\noindent
{\it MSC 2000.} Primary: 82C22; Secondary: 60K35, 92D25\\
%82C22    	Interacting particle systems 
%60K35    	Interacting random processes; statistical mechanics type
%               models; percolation theory
%92D25    	Population dynamics (general)
%
%60J80    	Branching processes (Galton-Watson, birth-and-death, etc.)
%60J60    	Diffusion processes
{\it Keywords.}
Reaction-diffusion process, branching, coalescence, annihilation, thinning,
Poissonization.\\
{\it Acknowledgments.} Work sponsored by GA\v{C}R grant: P201/10/0752.
Part of this work was done when the authors were visiting the Universit\"at
Duisburg-Essen. We would like to thank Anita Winter and the
staff for their hospitality;  Soumik Pal for
useful discussions; and an anonymous referee for a careful reading of the paper.
%First Schl\"ogl model, autocatalytic reaction,

{\setlength{\parskip}{-2pt}\tableofcontents}
\newpage

\section{Results}

\subsection{Introduction}\label{S:intro}

In \cite{AS05}, we studied systems of particles that perform independent
random walks, branch binarily, coalesce, and die. Our motivation came from two
directions. On the one hand, we were driven by the wish to study a population
dynamic model that is more realistic than the usual branching particle
systems, since the population at a given site cannot grow unboundedly but is
instead controlled by an extra death term that is quadratic in the number of
particles, which can be interpreted as extra deaths due to competition. On the
other hand, such systems of branching and coalescing particles are known to be
dual to certain systems of interacting diffusions, modelling gene frequencies
in spatially structured populations subject to resampling, mutation, and
selection \cite{SU86}. In this context, the branching-coalecing particles can
be interpreted as `potential ancestors' \cite{KN97}.

Apart from this duality, which was known, we showed in \cite{AS05} that our
particle systems are also related to resampling-selection processes by a
Poissonization relation. Moreover, we proved that systems started with
infinitely many particles on each site come down from infinity (a fact that
had been proved before, with a less explicit bound, in \cite{DDL90}) and that
systems on quite general spatially homogeneous lattices have at most one
nontrivial homogeneous invariant law, which, if it exists, is the long-time
limit law of the process started in any nontrivial homogeneous initial law.

In the present paper, we generalize all these results to systems where
moreover, with some positive rate, pairs of particles on the same site {\em
  annihilate} each other, resulting in the disappearance of {\em both}
particles. This my not seem like it should make a big difference with
coalescence, where only one particle disappears -and indeed our results
confirm this- but from the technical point of view annihilation has the huge
disadvantage of making the system non-monotone, which means that many simple
comparison arguments are not available. Some pioneering work on non-monotone
systems can be found in, e.g., \cite{BG85,Sud90,Dur91}. Despite progress in
recent years, non-monotone particle systems are still generally less studied
and worse understood than monotone ones.

As in the case without annihilation, our main tool is duality. In fact, it
turns out that systems with annihilation are dual to the same Markov process
(a system of interacting Wright-Fisher diffusions) as those without it, but
with a different (and more complicated) duality function. As a result, we
obtain Poissonization and thinning relations which show, among others, that
systems with annihilation can be obtained from systems without it by
independent thinning. We reported these duality and thinning relations
before (without proof) in \cite{Swa06a}.

The paper is organized a follows. In Section~\ref{S:defs} we define our model
and the dual system of interacting diffusions. In Section~\ref{S:DPT} we state
our duality result and show how this implies Poissonization and thinning
relations. Section~\ref{S:results} presents our main results, showing that the
system started with infinitely many particles comes down from infinity and
that systems started in a spatially homogeneous, nontrivial invariant law
converge to a unique homogeneous invariant law. Section~\ref{S:discus}
contains more discussion and an overview of our proofs, which are given in
Section~\ref{S:proofs}.

\subsection{Definition of the models}\label{S:defs}

Let $\La$ be a finite or countably infinite set and let $q(i,j)\geq 0$
$(i,j\in\La,\ i\neq j$) be the transition rates of a continuous time Markov
process on $\La$, the {\em underlying motion}, which jumps from site $i$ to
site $j$ with rate $q(i,j)$. For notational convenience, we set $q(i,i):=0$
$(i\in\La)$. We assume that the rates $q(i,j)$ are uniformly summable and (in
a weak sense) irreducible, and that the counting measure on $\La$ is an
invariant law for the underlying motion, i.e.:
\be\ba{rl}\label{qassum}
{\rm(i)}&\dis\sup_i\sum_jq(i,j)<\infty,\\[5pt]
{\rm(ii)}&\dis\forall\De\sub\La,\ \De\neq\emptyset,\La\ \;\exists 
i\in\De,\ j\in\La\beh\De\mbox{ such that }q(i,j)>0\mbox{ or }q(j,i)>0,\\[5pt]
{\rm(iii)}& \dis\sum_jq^\dgg(i,j)=\sum_jq(i,j)\ \forall i\in\La,
\mbox{ where }q^\dgg(i,j):=q(j,i).
\ec
Here and elsewhere sums and suprema over $i,j$ always run over $\La$,
unless stated otherwise.\med

\noi
{\em Branching-annihilating particle systems.} We now let $(\La,q)$ be as
above, fix rates $a,b,c,d\geq 0$, and consider systems of particles subject to
the following dynamics.
\bei
\item[$1^\circ$] Each particle jumps, independently of the others, from site
$i$ to site $j$ with rate $q(i,j)$.
\item[$2^\circ$] Each pair of particles, present on the same site,
annihilates with rate $2a$, resulting in the
disappearance of both particles. 
\item[$3^\circ$] Each particle splits with rate $b$ into two new
particles, created on the position of the old one. 
\item[$4^\circ$] Each pair of particles, present on the same site,
coalesces with rate $2c$, resulting in the
creation of one new particle on the position of the two old ones.
\item[$5^\circ$]
Each particle dies (disappears) with rate $d$. 
\eei
Let $X_t(i)$ denote the number of particles present at site $i\in\La$
and time $t\geq 0$. Then $X=(X_t)_{t\geq 0}$, with
$X_t=(X_t(i))_{i\in\La}$, is a Markov process with formal generator
\bc\label{Gdef}
Gf(x)&:=&\dis\sum_{ij}q(i,j)x(i)\{f(x+\de_j-\de_i)-f(x)\}
+a\sum_ix(i)(x(i)-1)\{f(x-2\de_i)-f(x)\}\\
&&\dis+b\sum_ix(i)\{f(x+\de_i)-f(x)\}
+c\sum_ix(i)(x(i)-1)\{f(x-\de_i)-f(x)\}\\
&&\dis+d\sum_ix(i)\{f(x-\de_i)-f(x)\},
\ec
where $\de_i(j):=1$ if $i=j$ and $\de_i(j):=0$ otherwise.  We call $X$ the
{\em $(q,a,b,c,d)$-branco-process}.

The process $X$ can be defined for finite initial states and also for
some infinite initial states in an appropriate Liggett-Spitzer space.
Following \cite{LS81}, we define
\be\label{LS}
\Ei_\ga(\La):=\{x\in\N^\La:\|x\|_\ga<\infty\},\quad
\mbox{with}\quad\|x\|_\ga:=\sum_i\ga_i|x(i)|,
\ee
where $\ga=(\ga_i)_{i\in\La}$ are strictly positive constants satisfying
\be\label{gacon}
\sum_i\ga_i<\infty\quand
\sum_j(q(i,j)+q^\dgg(i,j))\ga_j\leq K\ga_i\quad(i\in\La)
\ee
for some $K<\infty$. (Our assumptions on $q$ imply the existence of a $\ga$
satisfying (\ref{gacon}).)\med

\noi
{\em Resampling selection processes.} Let $(\La,q)$ be as before, let $r,s,m$
be nonnegative constants, and let $\Xc=(\Xc_t)_{t\geq 0}$ be the
$[0,1]^\La$-valued Markov process given by the unique pathwise solutions to
the infinite dimensional stochastic differential equation (SDE) (see
\cite{SU86,AS05}):
\bc\label{sde}
\di\Xc_t(i)&=&\dis\sum_jq(j,i)(\Xc_t(j)-\Xc_t(i))\,\di t
+s\Xc_t(i)(1-\Xc_t(i))\,\di t-m\Xc_t(i)\,\di t\\
&&\dis+\sqrt{2r\Xc_t(i)(1-\Xc_t(i))}\,\di B_t(i)
\qquad\qquad(t\geq 0,\ i\in\La),
\ec
where $(B(i))_{i\in\La}$ is a collection of independent Brownian motions.
The process $\Xc$ is a system of linearly interacting Wright-Fisher
diffusions, also known as stepping stone model, which can be used to model the
spatial distribution of gene frequencies in the presence of resampling,
selection, and mutation. Following \cite{AS05}, we call $\Xc$ the
resampling-selection process with underlying motion $(\La,q)$,
resampling rate $r$, selection rate $s$, and mutation rate $m$, or
shortly the {\em $(q,r,s,m)$-resem-process}.

\subsection{Duality, Poissonization, and thinning}\label{S:DPT}

We start with some notation. For $\phi,\psi\in[-\infty,\infty]^\La$, we write
\be
\li\phi,\psi\re:=\sum_i\phi(i)\psi(i)\qquad\mbox{and}
\qquad|\phi|:=\sum_i|\phi(i)|,
\ee
whenever the infinite sums are defined. For any $\phi:\La\to[-1,1]$ and
$x:\La\to\N$ we write
\be\label{phix}
\phi^x:=\prod_i\phi(i)^{x(i)}\quad\mbox{with}\quad 0^0:=1
\ee
whenever the infinite product converges and the limit does not depend
on the order of the coordinates. The following proposition
generalizes \cite[Theorem~1~(a)]{AS05}.
\bp{\bf(Duality)}\label{P:dual}
Assume that $a+c>0$ and let
\be\label{rsmabcd}
\al=a/(a+c),\quad r=a+c,\quad s=(1+\al)b,\quand m=\al b+d,
\ee
or equivalently
\be\label{abcdrsm}
a=\al r,\quad b=s/(1+\al),\quad c=(1-\al)r,\quand 
d=m-\al s/(1+\al).
\ee
Let $X$ be a $(q,a,b,c,d)$-branco-process with $X_0\in\Ei_\ga(\La)$ a.s.\ and
let $\Xc^\dgg$ be a $(q^\dgg,r,s,m)$-resem-process, independent of $X$.
Suppose that one or more of the following conditions are satisfied:
\be\label{aXX}
{\rm (i)}\ \al<1,\quad{\rm(ii)}\ |X_0|<\infty\ {\rm a.s.},
\quad{\rm(iii)}\ |\Xc^\dgg_0|<\infty\ {\rm a.s.}
\ee
Then
\be\label{dufo}
\E\big[(1-(1+\al)\Xc^\dgg_0)^{\txt X_t}\big]
=\E\big[(1-(1+\al)\Xc^\dgg_t)^{\txt X_0}\big]\qquad(t\geq 0),
\ee
where the infinite products inside the expectation are a.s.\ well-defined.
\ep

\noi
Proposition~\ref{P:dual}, together with a self-duality for
$(q,r,s,m)$-resem-processes described in \cite[Theorem~1~(b)]{AS05}, implies 
that $(q,a,b,c,d)$-branco-processes can be obtained as Poissonizations of 
re\-samp\-ling-selection processes, and as thinnings of each other, as we
explain now. (These thinning relations will prove useful several times in what
will follow. On the other hand, we have no application of the Poissonization
relations, but since they are very similar and closely related, we treat them
here as well.)

If $\phi$ is a $\half^\La$-valued random variable, then by definition
a Poisson measure with random intensity $\phi$ is an $\N^\La$-valued
random variable $\Pois(\phi)$ whose law is uniquely determined by
\be\label{Poisdef}
\E\big[(1-\psi)^{\txt\Pois(\phi)}\big]=\E\big[\ex{-\li\phi,\psi\re}\big]
\qquad(\psi\in[0,1]^\La),
\ee
where we allow for the case that $\ex{-\li\phi,\psi\re}=e^{-\infty}:=0$.
In particular, if $\phi$ is nonrandom, then the components
$(\Pois(\phi)(i))_{i\in\La}$ are independent Poisson distributed
random variables with intensity $\phi(i)$.

If $x$ and $\phi$ are random variables taking values in $\N^\La$ and
$[0,1]^\La$, respectively, then by definition a $\phi$-thinning of $x$
is an $\N^\La$-valued random variable $\Thin_\phi(x)$ whose law is
uniquely determined by
\be\label{Thindef}
\E\big[(1-\psi)^{\txt\Thin_\phi(x)}\big]=\E\big[(1-\phi\psi)^{\txt x}\big]
\qquad(\psi\in[0,1]^\La).
\ee
In particular, when $x$ and $\phi$ are nonrandom and $x=\sum_n\de_{i_n}$, then
a $\phi$-thinning of $x$ can be constructed as
$\Thin_\phi(x):=\sum_n\chi_n\de_{i_n}$ where the $\chi_n$ are independent
$\{0,1\}$-valued random variables with $\P[\chi_n=1]=\phi(i_n)$.  More
generally, if $x$ and $\phi$ are random, then we may construct $\Thin_\phi(x)$
in such a way that its conditional law given $x$ and $\phi$ is as in the
deterministic case. It is not hard to check that (\ref{Thindef}) holds more
generally for any $\psi\in[0,2]^\La$ provided $(1-\psi)^{\txt\Thin_\phi(x)}$
is a.s.\ well-defined.

\bp{\bf(Poissonization and thinning)}\label{P:Poisthin}
Fix $s,m\geq 0$, $r>0$, and $0\leq\bet\leq\al\leq 1$ such that
$m-\frac{\bet}{1+\bet}s\geq 0$. Let $X$ and $\ov X$ be the $(q,\al
r,\frac{1}{1+\al}s,(1-\al)r,m-\frac{\al}{1+\al}s)$-branco-process and
$(q,\bet r,\frac{1}{1+\bet}s,(1-\bet)r,m
-\frac{\bet}{1+\bet}s)$-branco-process,
respectively, and let $\Xc$ be the
$(q,r,s,m)$-resem-process. Then
\be\label{Pois}
\P\big[X_0\in\cdot\,\big]=\P\big[\Pois(\ffrac{s}{(1+\al)r}\Xc_0)\in\cdot\,\big]
\quad\mbox{implies}\quad
\P\big[X_t\in\cdot\,\big]=\P\big[\Pois(\ffrac{s}{(1+\al)r}\Xc_t)\in\cdot\,\big]
\qquad(t\geq 0).
\ee
and
\be\label{thin}
\P\big[X_0\in\cdot\,\big]
=\P\big[\Thin_{\frac{1+\bet}{1+\al}}(\ov X_0)\in\cdot\,\big]
\quad\mbox{implies}\quad
\P\big[X_t\in\cdot\,\big]
=\P\big[\Thin_{\frac{1+\bet}{1+\al}}(\ov X_t)\in\cdot\,\big]
\qquad(t\geq 0).
\ee
\ep
{\bf Proof} Formula~(\ref{Pois}) has been proved in case $\al=0$ in
\cite{AS05}. The general case can be derived along the same lines.
Alternatively, this can be derived from the case $\al=0$ using the
fact that
$\P[\Thin_{\frac{1}{1+\al}}(\Pois(\frac{s}{r}\Xc_t))\in\cdot\,]=
\P[\Pois(\frac{s}{(1+\al)r}\Xc_t)\in\cdot\,]$, and formula
(\ref{thin}), which we prove now.

If the initial laws of $X$ and $\ov X$ are related as in (\ref{thin})
and $\Xc^\dgg$ is a $(q^\dgg,r,s,m)$-resem-process started in
$\Xc_0=\phi$ with $|\phi|<\infty$, then by (\ref{dufo}),
\be\ba{l}\label{thincalc}
\dis\E\big[(1-(1+\al)\phi)^{\Thin_{\frac{1+\bet}{1+\al}}(\ov X_t)}\big]
=\E\big[(1-(1+\bet)\phi)^{\ov X_t}\big]
=\E\big[(1-(1+\bet)\Xc^\dgg_t)^{\ov X_0}\big]\\[5pt]
\dis=\E\big[(1-(1+\al)\Xc^\dgg_t)^{\Thin_{\frac{1+\bet}{1+\al}}(\ov X_0)}\big]
=\E\big[(1-(1+\al)\Xc^\dgg_t)^{X_0}\big]
=\E\big[(1-(1+\al)\phi)^{X_t}\big]\qquad(t\geq 0),
\ec
where we have used that by \cite[Lemma~20]{AS05} one has $|\Xc^\dgg_t|<\infty$
a.s.\ for each $t\geq 0$, which guarantees that the infinite products are
a.s.\ well-defined. Since (\ref{thincalc}) holds for all $\phi\in[0,1]^\La$
with $|\phi|<\infty$, (\ref{thin}) follows.\qed

\noi
As an immediate corollary of formula~(\ref{thin}), we have:
\bcor{\bf(Thinnings of processes without annihilation)}\label{C:allthin}
Let $a,b,c,d\geq 0$ and $a+c>0$. Let $X$ be the
$(q,a,b,c,d)$-branco-process, $\al:=\frac{a}{a+c}$, and let $\ov X$ be
the $(q,0,(1+\al)b,a+c,\al b+d)$-branco-process. Then
\be\label{thin2}
\P[X_0\in\cdot\,]=\P[\Thin_{\frac{1}{1+\al}}(\ov X_0)\in\cdot\,]
\quad\mbox{implies}\quad
\P[X_t\in\cdot\,]=\P[\Thin_{\frac{1}{1+\al}}(\ov X_t)\in\cdot\,]
\qquad(t\geq 0).
\ee
\ecor
In particular, each branco-process with a positive annihilation rate
can be obtained as a thinning of a process with zero annihilation
rate.

\subsection{Main results}\label{S:results}

Let $\ov\N=\N\cup\{\infty\}$ denote the one-point compactification of
$\N$, and equip $\ov\N^\La$ with the product topology. We say that
probability measures $\nu_n$ on $\ov\N^\La$ converge weakly to a limit
$\nu$, denoted as $\nu_n\Rightarrow\nu$, when $\int\nu_n(\di
x)f(x)\to\int\nu(\di x)f(x)$ for every $f\in\Ci(\ov\N^\La)$, the space
of continuous real functions on $\ov\N^\La$.

Our first main result shows that it is possible to start a
$(q,a,b,c,d)$-branco-process with infinitely many particles at each site.  We
call this the $(q,a,b,c,d)$-branco process started at infinity. This result
generalizes \cite[Theorem~2]{AS05}. For branching-coalescing particle systems
on $\Z^d$ with more general branching and coalescing mechanisms, but without
annihilation, a similar result has been proved in \cite{DDL90}.
\bt{\bf (The maximal process)}\label{T:max}
Assume that $a+c>0$. Then there exists an $\Ei_\ga(\La)$-valued process
$X^{(\infty)}=(X^{(\infty)}_t)_{t>0}$ with the following properties:
\med

\noi
{\bf (a)} For each $\eps>0$, $(X^{(\infty)}_t)_{t\geq\eps}$ is the
$(q,a,b,c,d)$-branco-process starting in $X^{(\infty)}_\eps$.\medskip

\noi
{\bf (b)} Set $r:=a+b+c-d$. Then
\be\label{explicit}
\E[X^{(\infty)}_t(i)]\leq\left\{\ba{cl}
\dis\frac{r}{(2a+c)(1-e^{-rt})}\quad&\mbox{if }r\neq 0,\\[5pt]
\dis\frac{1}{(2a+c)t}\quad&\mbox{if }r=0\ea\right.
\qquad(i\in\La).
\ee

\noi
{\bf (c)} If $X^{(n)}$ are $(q,a,b,c,d)$-branco-processes starting in
initial states $x^{(n)}\in\Ei_\ga(\La)$ such that
\be
x^{(n)}(i)\up\infty\quad\mbox{as}\quad n\up\infty\qquad(i\in\La),
\ee
then
\be\label{tomax}
\Li(X^{(n)}_t)\Asto{n}\Li(X^{(\infty)}_t)\qquad(t>0).
\ee
{\bf (d)} There exists an invariant measure $\ov\nu$ of the
$(q,a,b,c,d)$-branco-process such that
\be
\Li(X^{(\infty)}_t)\Asto{t}\ov\nu.
\ee
{\bf (e)} The measure $\ov\nu$ is uniquely characterised by
\be\label{extinct}
\int\ov\nu(\di x)(1-(1+\al)\phi)^x=\P^\phi[\exists t\geq 0\mbox{ such that }
\Xc^\dgg_t=0]\qquad(\phi\in[0,1]^\La,\ |\phi|<\infty),
\ee
where $\al:=a/(a+c)$ and $\Xc^\dgg$ denotes the $(q^\dgg,a+c,(1+\al)b,\al
b+d))$-resem-process started in $\phi$.\med

\noi
{\bf (f)} If $r,s,m,\al,\bet$ are as in Proposition~\ref{P:Poisthin} and
$X^{(\infty)}$ and $\ov X^{(\infty)}$ are the corresponding branco-processes
started at infinity, then
\be\label{thinmax}
\P[X^{(\infty)}_t\in\cdot\,]
=\P[\Thin_{\frac{1+\al}{1+\bet}}(\ov X^{(\infty)}_t)\in\cdot\,]
\qquad(t\geq 0).
\ee
A similar thinning relation holds between their long-time limit laws.
\et
If $a=0$, then it has been shown in \cite[Theorem~2~(e)]{AS05} that
$\ov\nu$ dominates any other invariant measure in the stochastic order, hence
$\ov\nu$ can righteously be called the {\em upper invariant measure} of the
process. In the general case, when we have annihilation, we do not know how to
compare $\ov\nu$ with other invariant measures in the stochastic order, and we
only work with the characterization of $\ov\nu$ in (\ref{extinct}).

To formulate our final result, we need some definitions. Let $(\La,q)$ be
our lattice with jump kernel of the underlying motion, as before. By definition,
an {\em automorphism} of $(\La,q)$ is a bijection $g:\La\to\La$ such that
$q(gi,gj)=q(i,j)$ for all $i,j\in\La$. We denote the group of all
automorphisms of $(\La,q)$ by ${\rm Aut}(\La,q)$. We say that a subgroup
$G\sub{\rm Aut}(\La,q)$ is {\em transitive} if for each $i,j\in\La$ there
exists a $g\in G$ such that $gi=j$. We say that $(\La,q)$ is {\em homogeneous}
if ${\rm Aut}(\La,q)$ is transitive. We define shift operators
$T_g:\N^\La\to\N^\La$ by
\be\label{shift}
T_gx(j):=x(g^{-1}j)\qquad(i\in\La,\ x\in\N^\La,\ g\in{\rm Aut}(\La,q)).
\ee
If $G$ is a subgroup of ${\rm Aut}(\La,q)$, then we say that a
probability measure $\nu$ on $\N^\La$ is {\em $G$-homogeneous} if
$\nu\circ T_g^{-1}=\nu$ for all $g\in G$. For example, if $\La=\Z^d$
and $q(i,j)=1_{\{|i-j|=1\}}$ (nearest-neighbor random walk), then the
group $G$ of translations $i\mapsto i+j$ ($j\in\La$) is a transitive
subgroup of ${\rm Aut}(\La,q)$ and the $G$-homogeneous probability
measures are the translation invariant probability measures.

The next theorem, which generalizes \cite[Theorem~4~(a)]{AS05}, is our main
result.
\bt{\bf (Convergence to the upper invariant measure)}\label{T:hom}
Assume that $(\La,q)$ is infinite and homogeneous, $G$ is a transitive
subgroup of ${\rm Aut}(\La,q)$, and $a+c>0$. Let $X$ be the
$(q,a,b,c,d)$-branco process started in a $G$-homogeneous nontrivial initial
law $\Li(X_0)$. Then $\Li(X_t)\Rightarrow\ov\nu$ as $t\to\infty$, where
$\ov\nu$ is the measure in (\ref{extinct}).
\et

\subsection{Discussion and outline}\label{S:discus}

The dualities in Proposition~\ref{P:dual} and \cite[Theorem~1~(b)]{AS05}, as
well as the Poissonization and thinning relations in
Proposition~\ref{P:Poisthin} play a central role in the present paper. These
relations, whose discovery was the starting point of the present work, are
similar to duality and thinning relations between general nearest-neighbor
interacting particle systems discovered by Lloyd and Sudbury in
\cite{SL95,SL97,Sud00}. In fact, as has been demonstrated in \cite[Prop.~6 and
  Lemma~7]{Swa06a} (see also the more detailed preprint of the same paper,
\cite[Prop~4.2 and Lemma~4.3]{Swa06b}), our relations can (at least formally)
be obtained as `local mean field' limits of (a special case of) the relations
of Lloyd and Sudbury. In \cite{SL97}, Lloyd and Sudbury observed that quite
generally, if two interacting particle systems have the same dual (whith a
special sort of duality relation as described in that article), then one is a
thinning of the other. This general principle is also responsible for the
Poissonization and thinning relations of our Proposition~\ref{P:Poisthin}.

The thinning relation in Corollary~\ref{C:allthin} is especially noteworthy,
since it allows us to compare non-monotone systems (which are generally hard
to study) with monotone systems. Also, the thinning relation (\ref{thinmax})
allows us to prove that the unique nontrivial homogeneous invariant measures
of $(q,\al r,\frac{1}{1+\al}s,(1-\al)r,m-\frac{\al}{1+\al}s)$-branco-processes
are monotone in $\al$ (w.r.t.\ to the stochastic order). Such sort of
comparison results between non-monotone systems are rarely available.  In
fact, these thinning relations suggest that the ergodic behavior of $(q,\al
r,\frac{1}{1+\al}s,(1-\al)r,m-\frac{\al}{1+\al}s)$-branco-processes (with
$r,s,m$ fixed but arbitrary $\al$) and the $(r,s,m)$-resem process should all
be `basically the same'.

It does not seem straightforward to make this claim rigorous,
however. The reason is that Poissonization or thinning can only produce
certain initial laws. Thus, an ergodic result for resampling-selection
processes, as has been proved in \cite{SU86}, only implies an ergodic result for
branching-annihilating particle systems started in initial laws that are
Poisson with random intensity, and likewise, the ergodic result for
branching-annihilating particle systems in \cite{AS05} implies our
Theorem~\ref{T:hom} only for special initial laws, that are thinnings of other
laws.

Our main tool for proving the statement for general initial laws is, like in
our previous paper, duality. In this respect, our methods differ from those in
\cite{DDL90}, which are based on entropy calculations, but are similar to
those used in, for example, \cite{SU86,BDD91,AS05,SS08}. The papers
\cite{SU86,AS05} are particularly close in spirit. The sort of cancellative
systems type duality that we have to use in the present paper is somewhat
harder to work with than the additive systems type duality in
\cite{SU86,AS05}. Earlier applications of this sort of `cancellative' duality
can be found in \cite{BDD91,SS08}.

The remainder of this paper is devoted to proofs. Proposition~\ref{P:dual} and
Theorems~\ref{T:max} and \ref{T:hom} are proved in Sections \ref{S:dual},
\ref{S:max} and \ref{S:hom}, respectively.

\section{Proofs}\label{S:proofs}

\subsection{Construction and approximation}

\subsection{Finite systems}

We denote the set of finite particle configurations by
$\Ni(\La):=\{x\in\N^\La:|x|<\infty\}$ and let
\be\label{Sidef}
\Si(\Ni(\La))
:=\{f:\Ni(\La)\to\R:|f(x)|\leq K|x|^k+M\mbox{ for some }K,M,k\geq 0\}
\ee
denote the space of real functions on $\Ni(\La)$ satisfying a polynomial
growth condition. Recall the definition of the operator $G$ from (\ref{Gdef}).
Generalizing \cite[Prop.~8]{AS05}, we have the following result. Below and in
what follows, we let $\P^x$ denote the law of the $(q,a,b,c,d)$-branco-process
started in $x$ and we let $\E^x$ denote expectation with respect to $\P^x$.
\bp{\bf(Finite branco-processes)}\label{P:finmart}
Let $X$ be the $(q,a,b,c,d)$-branco-process started in a finite state
$x$. Then $X$ does not explode. Moreover, with $z^{\li k\re}
:=z(z+1)\cdots(z+k-1)$, one has
\be\label{kmom}
\E^x\big[|X_t|^{\li k\re}\big]\leq |x|^{\li k\re}e^{kbt}
\qquad(k=1,2,\ldots,\ t\geq 0).
\ee
For each $f\in\Si(\Ni(\La))$, one has $Gf\in\Si(\Ni(\La))$ and $X$ solves the
martingale problem for the operator $G$ with domain $\Si(\Ni(\La))$.
\ep
{\bf Proof} The proof of \cite[Prop.~8]{AS05} carries over without a change.\qed

\detail{We prove here the fact, claimed in \cite{AS05}, that
$f\in\Si(\Ni(\La))$ implies $Gf\in\Si(\Ni(\La))$. Indeed, if
\[
|f(x)|\leq K|x|^k+M,
\]
then 
\[\ba{r@{\,}c@{\,}l}
|Gf(x)|&\leq&\dis\sum_{ij}q(i,j)x(i)|f(x+\de_j-\de_i)-f(x)|
+a\sum_ix(i)(x(i)-1)|f(x-2\de_i)-f(x)|\\
&&\dis+b\sum_ix(i)|f(x+\de_i)-f(x)|
+c\sum_ix(i)(x(i)-1)|f(x-\de_i)-f(x)|\\
&&\dis+d\sum_ix(i)|f(x-\de_i)-f(x)|,
\ea\]
where
\[
|f(x+\de_i)-f(x)|\leq|f(x+\de_i)|+|f(x)|\leq K(|x|+1)^k+M+K|x|^k+M
\leq K'|x|^k+M'
\]
for some $K',M'<\infty$. Treating other, similar factors in the same fashion
we obtain
\[\ba{r@{\,}c@{\,}l}
|Gf(x)|&\leq&\dis\big(\sup_i\sum_jq(i,j)\big)|x|\big(K'|x|^k+M'\big)
+a|x|(|x|-1)\big(K'|x|^k+M'\big)\\
&&\dis+b|x|\big(K'|x|^k+M'\big)
+c|x|(|x|-1)\big(K'|x|^k+M'\big)\\
&&\dis+d|x|\big(K'|x|^k+M'\big).
\ea\]}

\noi
We equip $\N^\La$ with the componentwise order, i.e., for two states $x,\ti
x\in\N^\La$, we write $x\leq\ti x$ if $x(i)\leq \ti x(i)$ for all $i\in\La$.
In \cite{AS05}, we made extensive use of monotonicity of branching-coalescing
particle systems. For systems with annihilation, most of these arguments do no
longer work. In fact, we can only prove the following fact.
\bl{\bf(Comparison of branco-processes)}\label{L:comp}
Let $X$ and $\ti X$ be the $(q,a,b,c,d)$-branco-process and the
$(q, 0,\ti b,\ti c,\ti d)$-branco-process started in finite initial states
$x$ and $\ti x$, respectively. Assume that
\be
x\leq\ti x,\quad b\leq\ti b,\quad a+c\geq\ti c,\quad d\geq\ti d.
\ee
Then $X$ and $\ti X$ can be coupled in such a way that
\be\label{comp}
X_t\leq\ti X_t\qquad(t\geq 0).
\ee
\el
{\bf Proof} This can be proved in the same way as \cite[Lemma~9]{AS05}, by
constructing a bivariate process $(B,W)$, say of black and white particles,
such that $X=B$ are the black particles and $\ti X=B+W$ are the black and
white particles together, with dynamics as described there, except that each
pair of black particles, present at the same site, is replaced with rate
$2(1-\tet)c$ by one black and one white particle, with rate $2(1-\tet)a$ by
two white particles, with rate $2\tet c$ by one black particle, and with rate
$2\tet a$ by one white particle, where $\tet:=\ti c/(a+c)$.\qed

\noi
We will often need to compare two $(q,a,b,c,d)$-branco-processes with the same
parameters but different initial states. A convenient way to do this is to use
coupling. Let $(Y^{01},Y^{11},Y^{10})$ be a trivariate process, in which
particles jump, die and give birth to particles of their own
type, and pairs of particles of the same type annihilate and coalesce in the
usual way of a $(q,a,b,c,d)$-branco-processes, and in addition, pairs of
particles of different types coalesce to one new particle with a type that
depends on its parents, according to the following rates:
\be\ba{l}\label{standcoup}
01+10\mapsto\hspace{7.5pt} 11\quad\mbox{at rate }r,\\[1pt]
01+11\mapsto\left\{\ba{@{}l}
10\quad\mbox{at rate }2a,\\[1pt]
11\quad\mbox{at rate }2c,\ea\right.\\[1pt]
\ec
and similarly $10+11\mapsto 01$ or $11$ at rate $2a$ resp.\ $2c$. Then it is
easy to see that, for any choice of the parameter $r\geq 0$, both
$X:=Y^{01}+Y^{11}$ and $X':=Y^{10}+Y^{11}$ are
$(q,a,b,c,d)$-branco-processes. We will call this the {\em standard coupling}
with parameter $r$. Note that if $a=0$, then $X_0\leq X'_0$ implies
$X_t\leq X'_t$ for all $t\geq 0$ but the same conclusion cannot be drawn if
$a>0$ because of the transition $01+11\mapsto 10$.

Let $X$ be the $(q,a,b,c,d)$-branco-process. It follows from
Proposition~\ref{P:finmart} that the semigroup $(S_t)_{t\geq 0}$ defined by
\be\label{Stdef}
S_tf(x):=\E^x[f(X_t)]\qquad(t\geq 0,\ x\in\Ni(\La),\ f\in\Si(\Ni(\La)))
\ee
maps $\Si(\Ni(\La))$ into itself. The semigroup gives first moments of
functions of our process. We will also need a covariance formula for functions
of our process, the general form of which is well-known. Below, for any
measure $\mu$ and function $f$, we write $\mu f:=\int f\di\mu$ whenever the
integral is well-defined, and we let $\cov_\mu(f,g):=\mu(fg)-(\mu f)(\mu g)$
denote the covariance of functions $f,g$ under $\mu$. Note that if $\mu$ is a
probability measure on $\Ni(\La)$, then $\mu S_tf=\int\mu(\di x)\E^x(f(X_t)]$,
i.e., $\mu S_t$ is the law at time $t$ of the $(q,a,b,c,d)$-branco-processes
started in the initial law $\mu$.

\bl{\bf(Covariance formula)}\label{L:cov}
Let $(S_t)_{t\geq 0}$ be the semigroup defined in (\ref{Stdef}) and let $\mu$
be a probability measure on $\Ni(\La)$ such that $\int\mu(\di x)|x|^k<\infty$
for all $k\geq 1$. Then, for each $t\geq 0$ and $f,g\in\Si(\Ni(\La))$, one has
\be\label{cov}
\cov_{\mu S_t}(f,g)=\cov_\mu(S_tf,S_tg)+2\int_0^t\mu S_{t-s}\Ga(S_sf,S_sg)\di s,
\ee
where $\Ga(f,g):=\ffrac{1}{2}\big(G(fg)-(Gf)g-f(Gg)\big)$ is given by
\bc\label{Gaform}
\dis 2\Ga(f,g)(x)
&=&\dis\sum_{ij}q(i,j)x(i)
\big(f(x+\de_j-\de_i)-f(x)\big)\big(g(x+\de_j-\de_i)-g(x)\big)\\
&&\dis+a\sum_ix(i)(x(i)-1)
\big(f(x-2\de_i)-f(x)\big)\big(g(x-2\de_i)-g(x)\big)\\
&&\dis+b\sum_ix(i)
\big(f(x+\de_i)-f(x)\big)\big(g(x+\de_i)-g(x)\big)\\
&&\dis+c\sum_ix(i)(x(i)-1)
\big(f(x-\de_i)-f(x)\big)\big(g(x-\de_i)-g(x)\big)\\
&&\dis+d\sum_ix(i)
\big(f(x-\de_i)-f(x)\big)\big(g(x-\de_i)-g(x)\big).
\ec
\el
{\bf Proof} Formula (\ref{cov}) is standard, but the details of the proof vary
depending on the Markov process under consideration. In the present case, we
can copy the proof of \cite[Prop.~2.2]{Swa09} almost without a change. We
start by noting that $fg\in\Si(\Ni(\La))$ for all $f,g\in\Si(\Ni(\La))$, hence
$\Ga(f,g):=\ffrac{1}{2}\big(G(fg)-(Gf)g-f(Gg)\big)$ is well-defined for all
$f,g\in\Si(\Ni(\La))$. It is a straightforward excercise to check that
$\Ga(f,g)$ is given by (\ref{Gaform}). Now (\ref{cov}) will follow from a
standard argument (such as given in \cite[Prop.~2.2]{Swa09}) provided we show
that
\bc
\dis\dif{s}S_s\big((S_tf)(S_ug)\big)&=&\dis S_sG\big((S_tf)(S_ug)\big),\\[5pt]
\dis\dif{t}S_s\big((S_tf)(S_ug)\big)&=&\dis S_s\big((GS_tf)(S_ug)\big),\\[5pt]
\dis\dif{u}S_s\big((S_tf)(S_ug)\big)&=&\dis S_s\big((S_tf)(GS_ug)\big)
\ec
for all $0\leq s,t,u$ and $f,g\in\Si(\Ni(\La))$. Let us say that a sequence of
functions $f_n\in\Si(\Ni(\La))$ converges `nicely' to a limit
$f\in\Si(\Ni(\La))$ if $f_n\to f$ pointwise and there exist constants
$K,M,k\geq 0$ such that $\sup_n|f_n(x)|\leq K|x|^k+M$. Then (\ref{kmom}) and
dominated convergence show that $f_n\to f$ `nicely' implies $S_tf_n\to S_tf$
`nicely'. Note also that if $f_n,f,g\in\Si(\Ni(\La))$ and $f_n\to f$ `nicely',
then $f_ng\to fg$ `nicely'. It is easy to check that
$Gf\in\Si(\Ni(\La))$ for all $f\in\Si(\Ni(\La))$. Since the
$(q,a,b,c,d)$-branco-process $X^x$ started in a deterministic initial state
$X^x_0=x\in\Ni(\La)$ solves the martingale problem for $G$, we have
\be
t^{-1}\big(S_tf(x)-f(x)\big)
=t^{-1}\int_0^t\E\big[Gf(X^x_s)\big]\di s
\underset{t\down 0}{\longrightarrow} Gf(x)
\qquad\big(x\in\Ni(\La)\big),
\ee
which by (\ref{kmom}) and the fact that $Gf\in\Si(\Ni(\La))$ implies that
$t^{-1}(S_tf-f)\to Gf$ `nicely' as $t\down 0$. Combining three facts, we
see that
\be\ba{l}
\dis\dif{s}S_s\big((S_tf)(S_ug)\big)
=\lim_{\eps\down 0}S_s(P_\eps-1)\big((S_tf)(S_ug)\big)
=S_sG\big((S_tf)(S_ug)\big),\\[5pt]
\dis\dif{t}S_s\big((S_tf)(S_ug)\big)
=\lim_{\eps\down 0}S_s\big(((P_\eps-1)S_tf)(S_ug)\big)
=S_s\big((GS_tf)(S_ug)\big),
\ec
and similarly for the derivative w.r.t.\ $u$, where we are using that if the
right-hand derivative of a continuous real function exists in each point and
depends continuously on $t$, then the function is continuously
differentiable (see, e.g., \cite[Excersise~17.24]{HS75}).\qed

\subsection{Infinite systems}\label{S:infsys}

Recall the definition of the Liggett-Spitzer space $\Ei_\ga(\La)$ from
(\ref{LS}). We let $\Ci_{\rm Lip}(\Ei_\ga(\La))$ denote the class of Lipschitz
functions on $\Ei_\ga(\La)$, i.e., $f:\Ei_\ga(\La)\to\R$ such that
$|f(x)-f(y)|\leq L\|x-y\|_\ga$ for some $L<\infty$.

The main result of this section is the following generalization of
\cite[Prop.~11]{AS05}.
\bp{\bf (Construction of branco-processes)}\label{P:constr}
Let $(S_t)_{t\geq 0}$ be the semigroup defined in (\ref{Stdef}). For each
$f\in\Ci_{\rm Lip}(\Ei_\ga(\La))$ and $t\geq 0$, the function $S_tf$ defined
in (\ref{Stdef}) can be extended to a unique Lipschitz function on
$\Ei_\ga(\La)$, also denoted by $S_tf$. There exists a unique (in
distribution) time-homogeneous Markov process with cadlag sample paths in the
space $\Ei_\ga(\La)$ equipped with the norm $\|\cdot\|_\ga$, such that
\be
\E^x[f(X_t)]=S_tf(x)\qquad(f\in\Ci_{\rm Lip}(\Ei_\ga(\La)),
\ x\in\Ei_\ga(\La),\ t\geq 0).
\ee
\ep

\noi
To prepare for the proof of Proposition~\ref{P:constr}, we start with the
following lemma, which generalizes \cite[Lemma~12]{AS05}.
\bl{\bf(Action of the semigroup on Lipschitz functions)}\label{L:Lip}
Let $(S_t)_{t\geq 0}$ be the semigroup of the $(q,a,b,c,d)$-branco-process,
defined in (\ref{Stdef}). If
$f:\Ni(\La)\to\R$ is Lipschitz continuous in the norm
$\|\cdot\|_\ga$ from (\ref{gacon}), with Lipschitz constant $L$, then
\be\label{Lipest}
|S_tf(x)-S_tf(x')|\leq Le^{(K+b-d)t}\|x-x'\|_\ga
\qquad(x,x'\in\Ni(\La),\ t\geq 0),
\ee
where $K$ is the constant from (\ref{gacon}).
\el
{\bf Proof} Let $X=Y^{01}+Y^{11}$ and $X'=Y^{10}+Y^{11}$ be
$(q,a,b,c,d)$-branco-processes started in $X_0=x$ and $X'_0=x'$, coupled using
the standard coupling from (\ref{standcoup}), in such a way that
$(Y^{01}_0,Y^{11}_0,Y^{10}_0)=((x-x')_+,x\wedge x',(x'-x)_+)$.
Then
\be\ba{l}\label{Lipa}
\dis|S_tf(x)-S_tf(x')|
=\big|\E[f(X_t)]-\E[f(X'_t)]\big|
\leq\E\big[|f(X_t)-f(X'_t)|\big]\\[5pt]
\dis\quad\leq L\E\big[\|X_t-X'_t\|_\ga\big]
=L\E\big[\|Y^{01}_t+Y^{10}_t\|_\ga\big].
\ec
Let us choose the parameter $r$ in the standard coupling as $r:=2(a+c)$. Then
it is easy to see that $(Y^{01},Y^{10})$ can be coupled to a
$(q,0,b,a+c,d)$-branco-process $Z$ started in $Z_0=|x-x'|$ in such a way that
$Y^{01}_t+Y^{10}_t\leq Z_t$ for all $t\geq 0$. Therefore, by
\cite[formula~(3.13)]{AS05}, we can further estimate the quantity in
(\ref{Lipa}) as
\be
|S_tf(x)-S_tf(x')|\leq L\E\big[\|Z_t\|_\ga\big]\leq Le^{(K+b-d)t}\|x-x'\|_\ga.
\ee
\qed

\noi
{\bf Proof of Proposition~\ref{P:constr}} Since $\Ni(\La)$ is a dense subset
of $\Ei_\ga(\La)$, Lemma~\ref{L:Lip} implies that for each $f\in\Ci_{\rm
  Lip}(\Ei_\ga(\La))$ and $t\geq 0$, the function $S_tf$ defined in
(\ref{Stdef}) can be extended to a unique Lipschitz function on
$\Ei_\ga(\La)$. The proof of Lemma~\ref{L:Lip} moreover shows that two
$(q,a,b,c,d)$-branco-processes $X,X'$ started in finite initial states $x,x'$
can be coupled such that
\be\label{xcoup}
\E\big[\|X_t-X'_t\|_\ga\big]\leq e^{(K+b-d)t}\|x-x'\|_\ga\qquad(t\geq 0).
\ee
It is not hard to see that for each $x\in\Ei_\ga(\La)$ we can choose
$x_n\in\Ni(\La)$ such that $\|x_n-x\|\to 0$ and
\be\label{sumxn}
\sum_{n=1}^\infty\|x_n-x_{n-1}\|_\ga<\infty.
\ee
(For example, any $x_n\up x$ has these properties.) Let $X^n$ be the process
started in $X^n_0=x_n$. By (\ref{xcoup}), we can inductively couple the
processes $X^0,X^1,X^2,\ldots$ in such a way that
\be
\E\big[\|X^n_t-X^{n-1}_t\|_\ga\big]\leq e^{(K+b-d)t}\|x_n-x_{n-1}\|_\ga
\qquad(n\geq 1,\ t\geq 0).
\ee
It follows that for each (deterministic) $t\geq 0$, the sequence
$(X^n_t)^{n\geq 0}$ is a.s.\ a Cauchy sequence in the complete metric space
$\Ei_\ga(\La)$, hence for each $t\geq 0$ there a.s.\ exists an
$\Ei_\ga(\La)$-valued random variable $X_t$ such that $\|X^n_t-X_t\|_\ga\to
0$. By Fatou,
\be
\E\big[\|X^n_t-X_t\|_\ga]\leq\liminf_{m\to\infty}\E\big[\|X^n_t-X^m_t\|_\ga]
\leq e^{(K+b-d)t}\sum_{m=n}^\infty\|x_{m+1}-x_m\|_\ga\asto{n}0.
\ee
Just as in \cite[Lemma~13]{AS05}, it is now straightforward to check that
$(X_t)_{t\geq 0}$ is a Markov process with semigroup $(S_t)_{t\geq 0}$. Note,
however, that in the arguments so far we have only constructed $X=(X_t)_{t\geq
  0}$ at deterministic times. To show that $X$ has a version with cadlag
sample paths (where the limits from the left and right are defined w.r.t.\ the
norm $\|\,\cdot\,\|_\ga$), we adapt arguments from the proof of
\cite[Prop.~11]{AS05}. It suffices to prove $X$ has cadlag sample paths on the
time interval $[0,1]$.

Let $V$ be the process with generator
\be\label{GV}
G_Vf(x):=\sum_{ij}q(i,j)x(i)\{f(x+\de_j)-f(x)\}+b\sum_ix(i)\{f(x+\de_i)-f(x)\},
\ee
which describes a branching process in which particles don't move or die, and
each particle at $i$ gives with rate $q(i,j)$ birth to a particle at $j$ and
with rate $b$ to a particle at $i$. We claim that a
$(q,a,b,c,d)$-branco-process $X$, started in a finite initial state $X_0=x$,
can be coupled to the process $V$ started in $V_0=x$ in such a way that
$X_t\leq V_t$ for all $t\geq 0$. To see this, let $(B,W)$ be a bivariate
process, say of black and white particles, started in $(B_0,W_0)=(x,0)$, such
that the black particles evolve as a $(q,a,b,c,d)$-branco-process, the white
particles evolve according to the generator in (\ref{GV}), and each time a
black particle disappears from a site $i$ due to jumps, annihilation or
coalescence, a white particle is created at~$i$. Then it is easy to see that
$X=B$ and $V=B+W$. By \cite[formula~(3.25)]{AS05},
\be
\E\big[\|V_t\|_\ga\big]\leq e^{(K+b)t}\|x\|_\ga,
\ee
where $K$ is the constant from (\ref{gacon}). Since $V$ is nondecreasing in
$t$, since $V_t(i)$ increases by one each time $X_t(i)$ does, and since $X$
cannot become negative, it follows that
\be\label{jumpbd}
\big|\big\{t\in[0,1]:X_{t-}(i)\neq X_t(i)\big\}\big|\leq x(i)+2V_1(i).
\ee
Applying this to the process $X^n$, multiplying with $\ga_i$ and summing over
$i$, we see that
\be\label{jmpbd}
\sum_i\ga_i\E\big[\big|\{t\in[0,1]:X^n_{t-}(i)\neq X^n_t(i)\}\big|\big]
\leq(1+2e^{K+b})\|x_n\|_\ga,
\ee
which by the convergence of $\|x_n\|_\ga$ gives us a uniform bound on the
number of jumps made by $X^n$.

We wish to show that for large $n$, the processes $X^n$ and $X^{n+1}$ make
mostly the same jumps. To this aim, let $X^n=Y^{01}+Y^{11}$ and
$X^{n+1}=Y^{10}+Y^{11}$ be two $(q,a,b,c,d)$-branco-processes, coupled using
the standard coupling from (\ref{standcoup}), with $r=2(a+c)$ and
$(Y^{01}_0,Y^{11}_0,Y^{10}_0)=((x_n-x_{n+1})_+,x_n\wedge
x_{n+1},(x_{n+1}-x)_+)$. Then, just as in the proof of Lemma~\ref{L:Lip}, the
process $(Y^{01},Y^{10})$ can be coupled to a $(q,0,b,a+c,d)$-branco-process
$Z$ started in $Z_0=|x_n-x_{n+1}|$ in such a way that $Y^{01}_t+Y^{10}_t\leq
Z_t$ for all $t\geq 0$. Likewise, it is not hard to see that we can couple
$(Y^{01},Y^{10})$ to a process $V$ with dynamics as in (\ref{GV}) started in
$V_0=|x_n-x_{n+1}|$, in such a way that $Y^{01}_t+Y^{10}_t\leq Z_t$ for all
$t\geq 0$ and moreover, whenever $Y^{01}(i)$ or $Y^{01}(i)$ increases, the
process $V(i)$ increases by the same amount. Let
\be
J_n(i):=\big\{t\in[0,1]:X^n_{t-}(i)\neq X^n_t(i)\big\}
\ee
be the set of jump times up to time one of the process $X^n(i)$ and let
\be
I(i):=\big\{t\in[0,1]:Y^{01}_{t-}(i)\neq Y^{01}_t(i)\big\}
\cup\big\{t\in[0,1]:Y^{10}_{t-}(i)\neq Y^{10}_t(i)\big\}.
\ee
Then the symmetric difference $J_n(i)\symdif J_{n+1}(i)=(J_n(i)\beh
J_{n+1}(i))\cup(J_{n+1}(i)\beh J_n(i))$ of $J_n(i)$ and $J_{n+1}(i)$ is
contained in $I(i)$ and, by the arguments leading up to (\ref{jumpbd}),
$|I(i)|\leq|x_n(i)-x_{n+1}(i)|+2V_1(i)$. Thus, in analogy with
(\ref{jmpbd}), we find that
\be\label{jmpdif}
\sum_i\ga_i\E\big[|J_n(i)\symdif J_{n+1}(i)|\big]
\leq(1+2e^{K+b})\|x_n-x_{n+1}\|_\ga.
\ee
By (\ref{sumxn}), it follows that the sets $J_n(i)$ converge as $n\to\infty$,
i.e., for each $i\in\La$ there is a (random) $n$ such that
$J_n(i)=J_{n+1}(i)=J_{n+2}(i)=\cdots$. Taking into account also
(\ref{jumpbd}), it follows that the limit process $(X(i))_{t\geq 0}$ has
cadlag sample paths for each $i\in\La$ and the set of jump times of $X^n(i)$
converges to the set of jump times of $X(i)$. The fact that the sample path of
$(X)_{t\geq 0}$ are also cadlag in the norm $\|\,\cdot\,\|_\ga$ can be proved
in the same way as \cite[formula~(3.31)]{AS05}.\qed

\noi
The proof of Proposition~\ref{P:constr} yields a useful side result.
\bcor{\bf(Approximation with finite systems)}\label{C:finapp}
Let $x\in\Ei_\ga(\La)$ and $x_n\in\Ni(\La)$ satisfy $\|x_n-x\|_\ga\to 0$ and
$\sum_{n\geq 1}\|x_n-x_{n-1}\|_\ga<\infty$. Then the
$(q,a,b,c,d)$-branco-processes $X^n,X$ started in $X^n_0=x_n$ and $X_0=x$ can
be coupled in such a way that $\|X^n_t-X_t\|_\ga\to 0$ a.s.\ for each $t\geq 0$.
\ecor

\subsection{Covariance estimates}

In this section, we give an upper estimate on the covariance of two functions
of a $(q,a,b,c,d)$-branco-process, which shows in particular that events that
are sufficiently far apart are almost independent.

For any continuous $f:\Ei_\ga(\La)\to\R$, we define $\de f:\La\to[0,\infty]$ by
\be
\de f(i):=\sup_{x\in\Ei_\ga(\La)}\big|f(x+\de_i)-f(x)\big|\qquad(i\in\La).
\ee
It is easy to see that for each continuous $f:\Ei_\ga(\La)\to\R$,
\be\label{fxyest}
\big|f(x)-f(y)\big|\leq\sum_i\de f(i)|x(i)-y(i)|
\qquad\big(x,y\in\Ei_\ga(\La)\big).
\ee

\bl{\bf(Lipschitz functions)}\label{L:Lipde}
A continuous function $f:\Ei_\ga(\La)\to\R$ is Lipschitz with respect to the
norm $\|\,\cdot\,\|_\ga$ if and only if there exists a constant $L<\infty$
such that $\de f(i)\leq L\ga_i$ $(i\in\La)$.
\el
{\bf Proof} If $f\in\Ci_{\rm Lip}(\Ei_\ga(\La))$, we have
$|f(x+\de_i)-f(x)|\leq L\|(x+\de_i)-x\|_\ga=L\ga_i$, where $L$ is the Lipschitz
constant of $f$, hence $\de f(i)\leq L\ga_i$ $(i\in\La)$. Conversely, if the
latter condition holds, then by (\ref{fxyest})
\be
\big|f(x)-f(y)\big|\leq L\sum_i\ga_i|x(i)-y(i)|=L\|x-y\|_\ga
\qquad\big(x,y\in\Ei_\ga(\La)\big).
\ee
\qed

We let $B_\ga(\La)$ denote the space of all functions $\phi:\La\to\R$ such
that
\be
\sup_i\ga_i^{-1}|\phi(i)|<\infty.
\ee
Note that by Lemma~\ref{L:Lipde}, $\de f\in B_\ga(\La)$ for each $f\in\Ci_{\rm
  Lip}(\Ei_\ga(\La))$.

Let $P_t(i,j)$ denote the probability that the random walk on $\La$ that jumps
from $k$ to $l$ with rate $q(k,l)$, started in $i$, is a time $t$ located at
the position $j$. For any $\phi\in B_\ga(\La)$, we write
\be
P_t\phi(i):=\sum_jP_t(i,j)\phi(j)\qquad(t\geq 0,\ i\in\La).
\ee
It is not hard to check that $P_t$ is well-defined on $B_\ga(\La)$ and maps
this space into itself.

Recall that $(S_t)_{t\geq 0}$ denotes the semigroup of the
$(q,a,b,c,d)$-branco-process, defined in (\ref{Stdef}).

\bl{\bf(Variation estimate)}\label{L:varest}
For any $(q,a,b,c,d)$-branco-process, one has
\be
\de S_tf\leq e^{(b-d)t}P_t\de f
\qquad\big(t\geq 0,\ f\in\Ci_{\rm Lip}(\Ei_\ga(\La))\big).
\ee
\el
{\bf Proof} Fix $i\in\La$ and let $X=Y^{01}+Y^{11}$ and $X'=Y^{10}+Y^{11}$ be
$(q,a,b,c,d)$-branco-processes started in $X_0=x$ and $X'_0=x+\de_i$, coupled
using the standard coupling from (\ref{standcoup}), in such a way that
$(Y^{01}_0,Y^{11}_0,Y^{10}_0)=(0,x,\de_i)$. Then
\be\ba{l}
\dis|S_tf(x)-S_tf(x+\de_i)|
=\big|\E[f(X_t)]-\E[f(X'_t)]\big|
\leq\E\big[|f(X_t)-f(X'_t)|\big]\\[5pt]
\dis\quad\leq\E\big[\sum_j\de f(j)|X_t(j)-X'_t(j)|\big]
=\sum_j\de f(j)\E\big[Y^{01}_t(j)+Y^{10}_t(j)\big]
\leq\sum_j\de f(j)e^{(b-d)t}P_t(i,j),
\ec
where in the last step we have used that $Y^{01}+Y^{10}$ can be estimated from
above by a $(q,0,b,0,d)$-branco-process.\qed

\bp{\bf(Covariance estimate)}\label{P:covest}
Let $X=(X_t)_{t\geq 0}$ be a $(q,a,b,c,d)$-branco-processes started in
$X_0=x\in\Ei_\ga(\La)$. Then, for each $t\geq 0$, there exist functions
$K_t:\La\times\La^2\to\half$ and $L_t:\La^2\times\La^2\to\half$ satisfying
\be\label{KLhom}
\left.\ba{l}
K_t(gi;gk,gl)=K_t(i;k,l)\\
L_t(gi,gj;gk,gl)=L_t(i,j;k,l)
\ea\right\}\quad\big(i,j,k,l\in\La,\ g\in{\rm Aut}(\La,q)\big),
\ee
\be\label{KLsum}
\mbox{and}\qquad\sup_{t\in[0,T]}\sum_{i,k}K_t(i;k,0)<\infty,
\quad\mbox{and}\quad
\sup_{t\in[0,T]}\sum_{i,j,k}L_t(i,j;k,0)<\infty\quad(T<\infty),
\ee
such that
\bc\label{covest}
\dis\big|\cov_x\big(f(X_t),g(X_t)\big)\big|
&\leq&\dis\sum_{i,k,l}x(i)K_t(i;k,l)\de f(k)\de g(l)\\[5pt]
&&\dis+\sum_{i,j,k.l}x(i)x(j)L_t(i,j;k,l)\de f(k)\de g(l).
\ec
for all bounded functions $f,g\in\Ci_{\rm Lip}(\Ei_\ga(\La))\big)$.
\ep
{\bf Proof} It suffices to prove the claim for finite initial states
$x\in\Ni(\La)$. For once the proposition is proved for finite systems, for
arbitrary $x\in\Ei_\ga(\La)$ we can find $\Ni(\La)\ni x_n\up x$. Then by
Corollary~\ref{C:finapp}, the processes $X^n,X$ started in $x_n,x$ can be
coupled such that $\|X^n_t-X_t\|_\ga\to 0$ for each $t\geq 0$, hence by
bounded pointwise convergence, the left-hand side of (\ref{covest}) for $X^n$
converges to the same formula for $X$, while the right-hand side is obviously
continuous under monotone limits.

We will show that for finite systems, the estimate (\ref{covest}) holds even
without the boundednes assumption on $f,g$. We apply Lemma~\ref{L:cov}. A
little calculation based on (\ref{Gaform}) shows that
\bc\label{Gaest}
\dis 2|\Ga(f,g)(x)|
&\leq&\dis\sum_{ij}q(i,j)x(i)
\big(\de f(i)+\de f(j)\big)\big(\de g(i)+\de g(j)\big)\\
&&\dis+(2a+c)\sum_ix(i)(x(i)-1)\de f(i)\de g(i)\\
&&\dis+(b+d)\sum_ix(i)\de f(i)\de g(i).
\ec
In view of Lemma~\ref{L:varest}, we define $\ti P_t:=e^{(b-d)t}P_t$.
Then (\ref{cov}), (\ref{Gaest}) and Lemma~\ref{L:varest} show that
for processes started in a deterministic initial state,
\be\ba{l}\label{fgX}
\dis\big|\cov_x\big(f(X_t),g(X_t)\big)\big|\\[5pt]
\dis\quad\leq\int_0^t\sum_{ij}q(i,j)
\big(\ti P_s\de f(i)+\ti P_s\de f(j)\big)
\big(\ti P_s\de g(i)+\ti P_s\de g(j)\big)
\E^x[X_{t-s}(i)]\di s\\[5pt]
\dis\quad\phantom{\leq}+(2a+c)\int_0^t\sum_i\ti P_s\de f(i)\ti P_s\de g(i)
\E^x\big[X_{t-s}(i)(X_{t-s}(i)-1)\big]\di s\\[5pt]
\dis\quad\phantom{\leq}+(b+d)\int_0^t\sum_i\ti P_s\de f(i)\ti P_s\de g(i)
\E^x[X_{t-s}(i)]\di s,
\ec
Let $Y=(Y_t)_{t\geq 0}$ be the $(q,0,b,0,d)$-branco-process started in
$Y_0=x$. By Lemma~\ref{L:comp}, we can couple $X$ and $Y$ such that $X_t\leq
Y_t$ for all $t\geq 0$. We estimate
\be\ba{rr@{\,}c@{\,}l}\label{XleqY}
{\rm(i)}&\dis\E^x[X_t(i)]&\leq&\dis\E^x[Y_t(i)]=\sum_jx(j)\ti P_t(j,i),\\[5pt]
{\rm(ii)}&\dis\E^x[X_t(i)(X_t(i)-1)]&\leq&\dis\E^x[Y_t(i)^2]
=\E^x[Y_t(i)]^2+\var_x(Y_t(i)).
\ec
To estimate $\var_x(Y_t(i))$, we apply (\ref{fgX}) to the process $Y$ and
$f=g:=f_i$ where $f_i(x):=x(i)$. Since the annihilation and coalescence rates
of $Y$ are zero, this yields
\be\ba{l}\label{fgY}
\dis\var_x(Y_t(i))\\[5pt]
\dis\quad\leq\int_0^t\sum_{jk}q(j,k)
\big(\ti P_s\de f_i(j)+\ti P_s\de f_i(k)\big)
\big(\ti P_s\de f_i(j)+\ti P_s\de f_i(k)\big)
\E^x[Y_{t-s}(j)]\di s\\[5pt]
\dis\quad\phantom{\leq}+(b+d)\int_0^t\sum_j\ti P_s\de f_i(j)\ti P_s\de f_i(j)
\E^x[Y_{t-s}(j)]\di s.
\ec
Define
\bc
\dis A_t(i;k,l)&:=&\dis\sum_jq(i,j)\big(\ti P_t(i,k)+\ti P_t(j,k)\big)
\big(\ti P_t(i,l)+\ti P_t(j,l)\big),\\[5pt]
\dis B_t(i;k,l)&:=&\dis\ti P_t(i,k)\ti P_t(i,l).
\ec
Then (\ref{fgX}) can be rewritten as
\be\ba{l}\label{Covrew}
\dis\big|\cov_x\big(f(X_t),g(X_t)\big)\big|
\leq\int_0^t\sum_{ikl}\E^x[X_{t-s}(i)]\big(A_s(i;k,l)+(b+d)B_s(i;k,l)\big)
\de f(k)\de_g(l)\di s\\[5pt]
\dis\qquad
+(2a+c)\int_0^t\sum_{ikl}\E^x\big[X_{t-s}(i)(X_{t-s}(i)-1)\big]B_s(i;k,l)
\de f(k)\de_g(l)\di s,
\ec
while (\ref{fgY}) can be rewritten as
\be\label{Varrew}
\var_x(Y_t(i))
\leq\int_0^t\sum_j\E^x[Y_{t-s}(j)]\big(A_s(j;i,i)+(b+d)B_s(j;i,i)\big)\di s,
\ee
where we have used that $\de f_i(j)=1_{\{i=j\}}$. Setting
\be
C_t(i;k,l):=A_s(i;k,l)+(b+d)B_s(i;k,l),
\ee
and inserting (\ref{XleqY}) and (\ref{Varrew}) into (\ref{Covrew}), we obtain
\be\ba{l}
\dis\big|\cov_x\big(f(X_t),g(X_t)\big)\big|\\[5pt]
\dis\leq\int_0^t\sum_{ijkl}x(i)\ti P_{t-s}(i,j)C_s(j;k,l)
\de f(k)\de_g(l)\di s\\[5pt]
\dis\phantom{\leq}+(2a+c)\int_0^t\sum_{ijklm}
x(i)x(j)\ti P_{t-s}(i,m)\ti P_{t-s}(j,m)B_s(m;k,l)
\de f(k)\de_g(l)\di s\\[5pt]
\dis\phantom{\leq}+(2a+c)\int_0^t\di s\sum_{ijklm}
\int_0^{t-s}\!\!\di u\, x(i)\ti P_u(i,j)
C_u(j;m,m)B_s(m;k,l)\de f(k)\de_g(l).
\ec
Recalling the definition of $B_t(i;j,k)$, this shows that (\ref{covest}) is
satisfied with
\bc\label{LKdef}
\dis K_t(i;k,l)
&:=&\dis\int_0^t\di s\sum_j\ti P_{t-s}(i,j)C_s(j;k,l)\\[5pt]
&&\dis+(2a+c)\int_0^t\di s\sum_{jm}
\int_0^{t-s}\!\!\di u\,\ti P_u(i,j)C_u(j;m,m)\ti P_s(m,k)\ti P_s(m,l),\\[5pt]
\dis L_t(i,j;k,l)
&:=&\dis(2a+c)\int_0^t\di s\sum_m
\ti P_{t-s}(i,m)\ti P_{t-s}(j,m)\ti P_s(m,k)\ti P_s(m,l).
\ec
The invariance of $K_t$ and $L_t$ under automorphisms of $(\La,q)$ is obvious
from the analogue property of $\ti P_t$, but the summability condition
(\ref{KLsum}) needs proof. Since $P_t(i,\,\cdot\,)$ is a probability
distribution and since the counting measure on $\La$ is an
invariant law for $P_t$ by assumption (\ref{qassum})~(iii),
\be
\sum_j P_t(i,j)=1=\sum_j P_t(j,i)\qquad(t\geq 0,\ i\in\La).
\ee
Setting $|q|:=\sum_jq(i,j)=\sum_jq(j,i)$, we see that
\bc
\dis\sum_{ik}A_t(i;k,l)
&=&\dis\sum_{ijk}q(i,j)\big(\ti P_t(i,k)+\ti P_t(j,k)\big)
\big(\ti P_t(i,l)+\ti P_t(j,l)\big)\\[5pt]
&=&\dis\sum_{ij}q(i,j)2e^{(b-d)t}
\big(\ti P_t(i,l)+\ti P_t(j,l)\big)=4|q|e^{2(b-d)t}\qquad(l\in\La),
\ec
and therefore, by a similar calculation for $B_t(i;k,l)$,
\be
\sum_jC_t(j;m,m)\leq
\sum_{jk}C_t(j;k,l)=\big(4|q|+b+d\big)e^{2(b-d)t}\qquad(l,m\in\La),
\ee
which by (\ref{LKdef}) implies that
\bc
\dis\sum_{ik}K_t(i;k,l)
&\leq&\dis\big(4|q|+b+d\big)\int_0^t\di s\, e^{(b-d)(t-s)}e^{2(b-d)s}\\[5pt]
&&\dis+(2a+c)\big(4|q|+b+d\big)\int_0^t\di s\int_0^{t-s}\!\!\di u\,
e^{(b-d)u}e^{2(b-d)u}e^{2(b-d)s}<\infty,\\[5pt]
\dis\sum_{ijk}L_t(i,j;k,l)
&\leq&\dis(2a+c)\int_0^t\di s\,e^{2(b-d)(t-s)}e^{2(b-d)s}<\infty
\qquad(t\geq 0,\ l\in\La).
\ec
\qed

\bcor{\bf(Exponential functions)}\label{C:exp}
Let $X=(X_t)_{t\geq 0}$ be a $(q,a,b,c,d)$-branco-processes started in
$X_0=x\in\Ei_\ga(\La)$, and let $\mu:\La\to\half$ satisfy $\sum_i\mu(i)<\infty$.
Then
\be\ba{l}\label{exp}
\dis\Big|\E^x\big[\ex{-\sum_i\mu(i)X_t(i)}\big]
-\prod_i\E^x\big[\ex{-\mu(i)X_t(i)}\big]\Big|\\[5pt]
\dis\quad\leq\ffrac{1}{2}\!\!\sum\subb{i,k,l}{k\neq l}\!
x(i)K_t(i;k,l)\mu(k)\mu(l)
+\ffrac{1}{2}\!\!\!\sum\subb{i,j,k,l}{k\neq l}\!\!
x(i)x(j)L_t(i,j;k,l)\mu(k)\mu(l),
\ec
where $K_t,L_t$ are as in Proposition~\ref{P:covest}.
\ecor
{\bf Proof} We first prove the statement if $\mu$ is finitely supported.
Let ${\rm support}(\mu)=\{k_1,\ldots,k_m\}$ and set
\be
f_\al(x):=\ex{-\mu(k_\al)x(k_\al)}
\quad\mbox{and}\quad
g_\bet(x):=\prod_{\al=1}^\bet f_\al(x).
\ee
Then
\bc
\dis\E^x\big[g_m(X_t)\big]
&=&\dis\E^x\big[g_{m-1}(X_t)\big]\E^x\big[f_m(X_t)\big]
+\cov_x\big(g_{m-1}(X_t),f_m(X_t)\big)\\[5pt]
&=&\dis\E^x\big[g_{m-2}(X_t)\big]\E^x\big[f_{m-1}(X_t)\big]
+\cov_x\big(g_{m-2}(X_t),f_{m-1}(X_t)\big)\\[3pt]
&&\dis\phantom{\E^x\big[g_{m-2}(X_t)\big]\E^x\big[f_{m-1}(X_t)\big]}
+\cov_x\big(g_{m-1}(X_t),f_m(X_t)\big)\\[5pt]
&=&\ \ldots\\[5pt]
&=&\dis\prod_{\al=1}^m\E[f_\al(X_t)]
+\sum_{\al=2}^m\cov_x\big(g_{\al-1}(X_t),f_\al(X_t)\big).
\ec
Therefore, since
\be
\de g_\al(k)=\left\{\ba{ll}
\mu(k)\quad&\mbox{if }k\in\{k_1,\ldots,k_\al\},\\
0\quad&\mbox{otherwise,}
\ea\right.\qquad
\de f_\al(k)=\left\{\ba{ll}
\mu(k)\quad&\mbox{if }k=k_\al,\\
0\quad&\mbox{otherwise,}\ea\right.
\ee
Proposition~\ref{P:covest} tells us that
\be\ba{l}
\dis\big|\E^x\big[g_m(X_t)\big]-\prod_{\al=1}^m\E[f_\al(X_t)]\big|\\[5pt]
\dis\quad\leq\sum_{\al=2}^m\sum_{\bet=1}^{\al-1}
\Big(\sum_ix(i)K_t(i;k_\bet,k_\al)\mu(k_\bet)\mu(k_\al)
+\sum_{i,j}x(i)x(j)L_t(i,j;k_\bet,k_\al)\mu(k_\bet)\mu(k_\al)\Big).
\ec
To generalize the statement to the case that $\sum_i\mu(i)<\infty$ but $\mu$
is not finitely supported, it suffices to choose finitely supported
$\mu_n\up\mu$ and to observe that all terms in (\ref{exp}) are continuous 
in $\mu$ w.r.t.\ increasing limits.\qed

\subsection{Duality and subduality}\label{S:dual}

Recall the definition of $\phi^x$ from (\ref{phix}).
\bl{\bf(Infinite products)}\label{L:infprod}
Let $0\leq\al\leq 1$, $\phi\in[0,1]^\La$, $x\in\N^\La$.\med

\noi
{\bf (a)} Assume that one or more of the following conditions are satisfied:
\be
{\rm (i)}\ \al<1,\quad{\rm(ii)}\ |x|<\infty,
\quad{\rm(iii)}\ |\phi|<\infty.
\ee
Then $\big(1-(1+\al)\phi\big)^x$ is well-defined.\med

\noi
{\bf (b)} Assume that $\phi$ is supported on a finite set and $x_n\in\N^\La$
converge pointwise to $x$. Then
$\big(1-(1+\al)\phi\big)^{x_n}\to\big(1-(1+\al)\phi\big)^x$ as $n\to\infty$.\med

\noi
{\bf (c)} Assume that $|\phi|<\infty$ and $\N^\La\ni x_n\up x$. Then
$\big(1-(1+\al)\phi\big)^{x_n}\to\big(1-(1+\al)\phi\big)^x$ as $n\to\infty$.\med

\noi
{\bf (d)} Assume that either $\al<1$ or $|\phi|<\infty$, and let
$[0,1]^\La\ni\phi_n\up\phi$. Then
$\big(1-(1+\al)\phi_n\big)^x\to\big(1-(1+\al)\phi\big)^x$ as $n\to\infty$.
\el
{\bf Proof} Since
$\big(1-(1+\al)\phi\big)^x:=\prod_i\big(1-(1+\al)\phi(i)\big)^{x(i)}$, where
$-1\leq\big(1-(1+\al)\phi(i)\big)^{x(i)}\leq 1$, the only way in which the
infinite product can be ill-defined is that
$\prod_i\big|1-(1+\al)\phi(i)\big|^{x(i)}>0$ while
$\big(1-(1+\al)\phi(i)\big)^{x(i)}<0$ for infinitely many $i$.
If $\al<1$, then $-1<-\al\leq 1-(1+\al)\phi(i)$, so if
$\big(1-(1+\al)\phi(i)\big)^{x(i)}<0$ for infinitely many $i$, then
$\prod_i\big|1-(1+\al)\phi(i)\big|^{x(i)}=0$ and the infinite product is
always well-defined. If $|x|<\infty$, then
$\big(1-(1+\al)\phi(i)\big)^{x(i)}=1$ for all but finitely many $i$, hence the
infinite product is certainly well-defined. If $|\phi|<\infty$, finally, then
$\phi(i)>\frac{1}{2}$ for finitely many $i$, hence
$\big(1-(1+\al)\phi(i)\big)^{x(i)}<0$ for finitely many $i$ and the infinite
product is again well-defined. This completes the proof of part~(a).

Part~(b) is trivial since all but finitely many factors in the infinite
product defining $\big(1-(1+\al)\phi\big)^x$ are one if $\phi$ is finitely
supported.

To prove part~(c), we split the product
$\prod_i\big(1-(1+\al)\phi(i)\big)^{x(i)}$ in finitely many factors where
$\phi(i)>\ffrac{1}{2}$ and the remaining factors where
$\phi(i)\leq\ffrac{1}{2}$ and hence $\big(1-(1+\al)\phi(i)\big)\geq 0$. Then
the finite part of the product converges as in part~(b) while the infinite
part converges in a monotone way.

For the proof of part~(d) set $I:=\{i\in\La:x(i)\neq 0\}$ and let
$I_-,I_0,I_+$ be the subsets of $I$ where $1-(1+\al)\phi(i)<0$, $=0$ and $>0$,
respectively. If $I_0\neq\emptyset$ then it is easy to see that 
$\big(1-(1+\al)\phi_n\big)^x\to 0=\big(1-(1+\al)\phi\big)^x$, so from now on
we may assume that $I_0=\emptyset$. Note that $1-(1+\al)\phi_n(i)\geq
1-(1+\al)\phi(i)>0$ for all $i\in I_+$. Therefore, if $I_-$ is finite,
as must be the case when $|\phi|<\infty$, then $\prod_{i\in
  I_-}\big(1-(1+\al)\phi_n(i)\big)^{x(i)}$ converges since $I_-$ is finite
while $\prod_{i\in I_+}\big(1-(1+\al)\phi_n(i)\big)^{x(i)}\down\prod_{i\in
  I_+}\big(1-(1+\al)\phi(i)\big)^{x(i)}$. If $I_-$ is infinite and $\al<1$,
then the fact that $\big|1-(1+\al)\phi_n(i)\big|^{x(i)}\to
\big((1+\al)\phi(i)-1\big)^{x(i)}\leq\al$ for each $i\in I_-$ implies that
$\big(1-(1+\al)\phi_n\big)^x\to 0=\big(1-(1+\al)\phi\big)^x$.\qed

\noi
We equip the space $[0,1]^\La$ with the product topology and let
$\Ci([0,1]^\La)$ denote the space of continuous real functions on
$[0,1]^\La$, equipped with the supremum norm. By $\Ci^2_{\rm
fin}([0,1]^\La)$ we denote the space of $\Ci^2$ functions on $[0,1]^\La$
depending on finitely many coordinates. By definition, $\Ci^2_{\rm
sum}([0,1]^\La)$ is the space of continuous functions $f$ on
$[0,1]^\La$ such that the partial derivatives $\dif{\phi(i)}f(\phi)$ and
$\difif{\phi(i)}{\phi(j)}f(\phi)$ exist for each $x\in(0,1)^\La$ and such that
the functions
\be\ba{l}
\phi\mapsto\big(\dif{\phi(i)}f(\phi)\big)_{i\in\La}
\quand
\phi\mapsto\big(\difif{\phi(i)}{\phi(j)}f(\phi)\big)_{i,j\in\La}
\ec
can be extended to continuous functions from $[0,1]^\La$ into the
spaces $\ell^1(\La)$ and $\ell^1(\La^2)$ of absolutely summable
sequences on $\La$ and $\La^2$, respectively, equipped with the
$\ell^1$-norm. Define an operator $\Gi:\Ci^2_{\rm
sum}([0,1]^\La)\to\Ci([0,1]^\La)$ by
\bc\label{Gidef}
\Gi f(\phi)&:=&\dis\sum_{ij}q(j,i)(\phi(j)-\phi(i))\dif{\phi(i)}f(\phi)
+s\sum_i\phi(i)(1-\phi(i))\dif{\phi(i)}f(\phi)\\[2pt]
&&\dis+r\sum_i\phi(i)(1-\phi(i))\diff{\phi(i)}f(\phi)
-m\sum_i\phi(i)\dif{\phi(i)}f(\phi)\qquad(\phi\in[0,1]^\La).
\ec
One can check that for $f\in\Ci^2_{\rm sum}([0,1]^\La)$, the infinite sums in
(\ref{Gidef}) converge in the supremumnorm and the result does not depend on
the summation order \cite[Lemma~3.4.4]{Swa99}. It has been shown in
\cite[Section~3.4]{AS05} that solutions to the SDE (\ref{sde}) solve the
martingale problem for the operator $\Gi$. In view of this, we loosely refer
to $\Gi$ as the {\em generator} of the $(q,r,s,m)$-resem-process.\med

\noi
{\bf Proof of Proposition~\ref{P:dual}} Since by Proposition~\ref{P:finmart}
(resp\ \cite[Lemma~20]{AS05}), $|X_0|<\infty$ (resp.\ $|\Xc_0|<\infty$) implies
$|X_t|<\infty$ (resp.\ $|\Xc_t|<\infty$) for all $t\geq 0$, by
Lemma~\ref{L:infprod}, each of the conditions (\ref{aXX})~(i)--(iii)
guarantees that both sides of equation (\ref{dufo}) are well-defined.

It suffices to prove (\ref{dufo}) for deterministic initial states, i.e., we
want to prove that either $\al<1$, $|x|<\infty$, or $|\phi|<\infty$ imply that
\be\label{dufo2}
\E^{\txt x}\big[(1-(1+\al)\phi)^{\txt X_t}\big]
=\E^{\txt\phi}\big[(1-(1+\al)\Xc^\dgg_t)^{\txt x}\big]\qquad(t\geq 0),
\ee
where $\E^x$ and $\E^\phi$ denote expectation w.r.t.\ the law of the process
$X$ started in $X_0=x$ and the process $\Xc$ started in $\Xc_0=\phi$,
respectively. We start by proving (\ref{dufo2}) if $|x|<\infty$. We wish to
apply \cite[Thm~7]{AS05}. Unfortunately, the original formulation of this
theorem contains an error, so we have to use the corrected version in
\cite[Corollary~2]{AS09b} (see also \cite[Corollary~2]{AS09a}). We apply this
to the duality function
\be
\Psi(x,\phi):=\big(1-(1+\al)\phi\big)^x
\qquad\big(x\in\Ni(\La),\ \phi\in[0,1]^\La\big).
\ee
Since $|\Psi(x,\phi)|\leq 1$, we obviously have
$\Psi(\,\cdot\,,\phi)\in\Si(\Ni(\La))$ for each $\phi\in[0,1]^\La$. Since for
each $x\in\Ni(\La)$, the function $\Psi(x,\,\cdot\,)$ depends only on finitely
many coordinates, we moreover have $\Psi(x,\,\cdot\,)\in\Ci^2_{\rm
  sum}([0,1]^\La)$ for each such $x$. Let $G$ be the generator of the
$(q,a,b,c,d)$-branco-process and let $\Gi^\dgg$
denote the generator of the $(q^\dgg,r,s,m)$-resem-process. In order to apply
\cite[Corollary~2]{AS09b}, we need to check that
\be\label{ducheck}
\Phi(x,\phi):=G\Psi(\,\cdot\,,\phi)(x)=\Gi^\dgg\Psi(x,\,\cdot\,)(\phi)
\qquad\big(x\in\Ni(\La),\ \phi\in[0,1]^\La\big)
\ee
and
\be\label{Edu}
\int_0^T\!\di s\,\int_0^T\!\di t\,
\E\big[\big|\Phi(X_s,\Xc_t)\big|\big]<\infty
\qquad(T\geq 0).
\ee
To check (\ref{ducheck}), we calculate:
\bc\label{check1}
G\Psi(\,\cdot\,,\phi)(x)&=&
\dis\sum_{ij}q(i,j)x(i)(1-(1+\al)\phi)^{x-\de_i}((1-(1+\al)\phi)^{\de_j}-(1-(1+\al)\phi)^{\de_i})\\[5pt]
&&\dis+a\sum_ix(i)(x(i)-1)(1-(1+\al)\phi)^{x-2\de_i}(1-(1-(1+\al)\phi)^{2\de_i})\\[5pt]
&&\dis+b\sum_ix(i)(1-(1+\al)\phi)^{x-\de_i}\big((1-(1+\al)\phi)^{2\de_i}-(1-(1+\al)\phi)^{\de_i}\big)\\[5pt]
&&\dis+c\sum_ix(i)(x(i)-1)(1-(1+\al)\phi)^{x-2\de_i}\big((1-(1+\al)\phi)^{\de_i}-(1-(1+\al)\phi)^{2\de_i}\big)\\[5pt]
&&\dis+d\sum_ix(i)(1-(1+\al)\phi)^{x-\de_i}(1-(1-(1+\al)\phi)^{\de_i}).
\ec
Since
\bc
\dis\dif{\phi(i)}(1-(1+\al)\phi)^x
&=&\dis-(1+\al)x(i)(1-(1+\al)\phi)^{x-\de_i},\\[5pt]
\dis\diff{\phi(i)}(1-(1+\al)\phi)^x
&=&\dis(1+\al)^2x(i)(x(i)-1)(1-(1+\al)\phi)^{x-2\de_i}
\ec
and
\bc
\dis(1-(1+\al)\phi)^{\de_i}&=&\dis 1-(1+\al)\phi(i),\\[5pt]
\dis(1-(1+\al)\phi)^{2\de_i}&=&\dis\big(1-(1+\al)\phi(i)\big)^2,\\[5pt]
\dis(1-(1+\al)\phi)^{\de_i}-(1-(1+\al)\phi)^{2\de_i}
&=&\dis(1+\al)\phi(i)\big(1-(1+\al)\phi(i)\big),
\ec
we can rewrite the expression in (\ref{check1}) as
\bc
G\Psi(\,\cdot\,,\phi)(x)&=&\dis\sum_{ij}q(i,j)(\phi(j)-\phi(i))\dif{\phi(i)}(1-(1+\al)\phi)^x\\[5pt]
&&\dis+\frac{a}{(1+\al)^2}\big(2(1+\al)\phi(i)-(1+\al)^2\phi(i)^2\big)\diff{\phi(i)}(1-(1+\al)\phi)^x\\[5pt]
&&\dis+\frac{b}{1+\al}(1+\al)\phi(i)\big(1-(1+\al)\phi(i)\big)\dif{\phi(i)}(1-(1+\al)\phi)^x\\[5pt]
&&\dis+\frac{c}{(1+\al)^2}(1+\al)\phi(i)\big(1-(1+\al)\phi(i)\big)\diff{\phi(i)}(1-(1+\al)\phi)^x\\[5pt]
&&\dis-\frac{d}{1+\al}(1+\al)\phi(i)\dif{\phi(i)}(1-(1+\al)\phi)^x.
\ec
Reordering terms gives
\bc
G\Psi(\,\cdot\,,\phi)(x)&=&\dis\sum_{ij}q(i,j)(\phi(j)-\phi(i))\dif{\phi(i)}(1-(1+\al)\phi)^x\\[5pt]
&&\dis+\Big(\frac{2a+c}{1+\al}\phi(i)-(a+c)\phi(i)^2\Big)\diff{\phi(i)}(1-(1+\al)\phi)^x\\[5pt]
&&\dis+\Big((b-d)\phi(i)-b(1+\al)\phi(i)^2\Big)\dif{\phi(i)}(1-(1+\al)\phi)^x\\[10pt]
&=&\dis\dis\sum_{ij}q^\dgg(j,i)(\phi(j)-\phi(i))\dif{\phi(i)}(1-(1+\al)\phi)^x\\[5pt]
&&\dis+(a+c)\phi(i)(1-\phi(i))\diff{\phi(i)}(1-(1+\al)\phi)^x\\[5pt]
&&\dis+(1+\al)b\phi(i)(1-\phi(i))\dif{\phi(i)}(1-(1+\al)\phi)^x\\[5pt]
&&\dis-(\al b+d)\phi(i)\dif{\phi(i)}(1-(1+\al)\phi)^x
\quad\qquad=\ \Gi^\dgg\Psi(x,\,\cdot\,)(\phi),
\ec
where we have used (\ref{rsmabcd}), which implies in particular that
\be
\frac{2a+c}{1+\al}=\frac{2a+c}{1+a/(a+c)}=\frac{(2a+c)(a+c)}{(a+c)+a}=a+c.
\ee
It is easy to see from (\ref{check1}) that there exists a constant $K$ such
that
\be
|\Phi(x,\phi)|\leq K\big(1+|x|^2\big)
\qquad\big(\phi\in[0,1]^\La,\ x\in\Ni(\La)\big),
\ee
hence (\ref{Edu}) follows from Proposition~\ref{P:finmart}. This completes the
proof of (\ref{dufo2}) in case $|x|<\infty$.

We next claim that (\ref{dufo2}) holds if $x\in\Ei_\ga(\La)$ and $\phi$ is
supported on a finite set. Choose $\Ni(\La)\ni x_n\up x$ and let $X^n$ denote
the $(q,a,b,c,d)$-branco-process started in $X^n_0=x_n$. Then
Corollary~\ref{C:finapp} implies that the $X^n$ can be coupled such that
$X^n_t(i)\to X_t(i)$ a.s.\ for each $i\in\La$. Therefore, taking the limit in
(\ref{dufo2}), using the fact that the integrands on the left- and right-hand
sides converge in a bounded pointwise way by Lemma~\ref{L:infprod}~(b) and
(c), respectively, our claim follows.

To see that (\ref{dufo2}) holds more generally if $\al<1$ or $|\phi|<\infty$,
we choose finitely supported $\phi_n\up\phi$ and let $\Xc^n$ denote the
$(q,r,s,m)$-resem-process started in $\Xc^n_0=\phi_n$. Then
\cite[Lemma~22]{AS05} implies that the $\Xc^n$ can be coupled such that
$\Xc^n_t(i)\up\Xc_t(i)$ a.s.\ for each $i\in\La$. The statement then follows
by letting $n\to\infty$ and applying Lemma~\ref{L:infprod}~(d).\qed

\noi
Fix constants $\bet\in\R$, $\ga\geq 0$. Let
$\Mi(\La):=\{\phi\in\half^\La:|\phi|<\infty\}$ be the space of finite
measures on $\La$, equipped with the topology of weak convergence, and
let $\Yi$ be the Markov process in $\Mi(\La)$ given by the unique
pathwise solutions to the SDE
\be\label{supersde}
\di\Yi_t(i)=\sum_ja(j,i)(\Yi_t(j)-\Yi_t(i))\,\di t
+\bet\Yi_t(i)\,\di t+\sqrt{2\ga\Yi_t(i)}\,\di B_t(i)
\ee
$(t\geq 0,\ i\in\La)$. Then $\Yi$ is the well-known super random walk
with underlying motion $a$, growth parameter
$\bet$ and activity $\ga$. One has \cite[Section~4.2]{Daw93}
\be\label{laplace}
\E^\phi\big[\ex{-\li\Yi_t,\psi\re}]=\ex{-\li\phi,\Ui_t\psi\re}
\ee
for any $\phi\in\Mi(\La)$ and bounded nonnegative $\psi:\La\to\R$, where
$u_t=\Ui_t\psi$ solves the semilinear Cauchy problem
\be\label{udif}
\dif{t}u_t(i)=\sum_ja(j,i)(u_t(j)-u_t(i))+\bet u_t(i)
-\ga u_t(i)^2\qquad(i\in\La,\ t\geq 0)
\ee
with initial condition $u_0=\psi$. The semigroup $(\Ui_t)_{t\geq 0}$
acting on bounded nonnegative functions $\psi$ on $\La$ is called the
log-Laplace semigroup of $\Yi$.

It has been shown in \cite[Prop.~23]{AS05} that the
$(q,a,b,c,d)$-branco-process and the super random walk with underlying motion
$q^\dgg$, growth parameter $b-d+c$ and activity $c$ are related by a
`subduality', i.e., a duality formula with a nonnegative error term.
The next proposition generalizes this to branco-processes with positive
annihilation rate.
\bp{\bf(Subduality with a branching process)}\label{P:subdup}
Let $X$ be the $(q,a,b,c,d)$-branco-process and let $\Yi$ be the super
random walk with underlying motion $q^\dgg$, growth parameter $2a+b-d+c$
and activity $2a+c$. Then
\be\label{subdu}
\E^x\big[\ex{-\li\phi,X_t\re}]\geq \E^\phi\big[\ex{-\li\Yi_t,x\re}]
\qquad(x\in\Ei_\ga(\La),\ \phi\in\half^\La,\ |\phi|<\infty).
\ee
\ep
{\bf Proof} We first prove the statement if $|x|<\infty$ and
$|\phi|<\infty$. This goes exactly in the same way as in the proof of
\cite[Prop.~23]{AS05}. Let $\Hi$ denote the generator of $\Yi$, defined in
\cite[formula~(4.14)]{AS05}, let $G$ be the generator in (\ref{Gdef}), and let
$\Psi$ be the duality function $\Psi(x,\phi):=e^{-\li\phi,x\re}$.
Then one has
\be\ba{l}
\dis G\Psi(\cdot,\phi)(x)-\Hi\Psi(x,\cdot)(\phi)=\Big\{\sum_{ij}q(i,j)x(i)
\big(e^{\phi(i)-\phi(j)}-1-(\phi(i)-\phi(j))\big)\\[5pt]
\quad\dis+a\sum_ix(i)(x(i)-1)\big(e^{2\phi(i)}-1-2\phi(i)\big)
+b\sum_ix(i)\big(e^{-\phi(i)}-1+\phi(i)\big)\\[5pt]
\quad\dis+c\sum_ix(i)(x(i)-1)\big(e^{\phi(i)}-1-\phi(i)\big)
+d\sum_ix(i)\big(e^{\phi(i)}-1-\phi(i)\big)\Big\}
\ex{-\li\phi,x\re}\geq 0.
\ec
This is just \cite[formula~(4.19)]{AS05}, where the extra terms with the
prefactor $a$ obtain their $e^{2\phi(i)}-1$ part from the generator $G$ and
the remaining $-2\phi(i)$ from $\Hi$. Using Proposition~\ref{P:finmart} to
guarantee integrability we may apply \cite[Corollary~2]{AS09b} to deduce
(\ref{subdu}).

To generalize (\ref{subdu}) to $x\in\Ei_\ga(\La)$ and $\phi\in\half^\La$
supported on a finite set, we choose $\Ni(\La)\ni x_n\up x$ and let $X^n$
denote the $(q,a,b,c,d)$-branco-process started in $X^n_0=x_n$. Then
Corollary~\ref{C:finapp} implies that the $X^n$ can be coupled such that
$X^n_t(i)\to X_t(i)$ a.s.\ for each $i\in\La$. It follows that
$e^{-\li\phi,X^n_t\re}\to e^{-\li\phi,X_t\re}$ a.s.\ and
$e^{-\li\Yi_t,x_n\re}\down e^{-\li\Yi_t,x_n\re}$ a.s., so taking the limit in
(\ref{subdu}) we obtain the statement for $x\in\Ei_\ga(\La)$ and $\phi$
finitely supported. To generalize this to $|\phi|<\infty$ we choose
$\phi_n\up\phi$ and let $\Yi^n$ denote the super random walk started in
$\Yi^n_0=\phi_n$. Then it is well-known (and can be proved in the same way as
\cite[Lemma~22]{AS05}) that the $\Yi^n$ can be coupled in such a way that
$\Yi^n_t\up\Yi_t$ for each $t\geq 0$. Therefore, taking the monotone
limit in (\ref{subdu}) our claim follows.\qed

\subsection{The process started at infinity}\label{S:max}

In view of what follows, we recall the following projective limit theorem.
Let $E$ and $(E_i)_{i\in\N}$ be Polish spaces. Assume that $\pi_i:E\to E_i$
are continuous surjective maps that separate points, i.e., for all $x,y\in E$
with $x\neq y$, there exists an $i\in\N$ with $\pi_i(x)\neq\pi_i(y)$. For each
$i\leq j$, let $\pi_{ij}:E_j\to E_i$ be continuous maps satisfying
$\pi_{ij}\circ\pi_j=\pi_i$. Assume moreover that for each sequence
$(x_i)_{i\in\N}$ with $x_i\in E_i$ $(i\in\N)$ that is consistent in the sense
that $\pi_{ij}(x_j)=x_i$ $(i\leq j)$, there exists an $y\in E$ such that
$\pi_i(y)=x_i$ $(i\in\N)$. Let $(\mu_i)_{i\in\N}$ be probability measures on
the $E_i$'s, respectively (equipped with their Borel-\si-fields), that are
consistent in the sense that $\mu_i=\mu_j\circ\pi_{ij}^{-1}$ for all $i\leq
j$. Then there exists a unique probability measure $\mu$ on $E$ such that
$\mu\circ\pi_i^{-1}=\mu_i$ for all $i\in\N$.

This may be proved by invoking Kolmogorov's extension theorem to construct a
probability measure $\mu'$ on the product space $\prod_iE_i$ whose marginals
are the $\mu_i$ and that is moreover concentrated on the set
$E'\sub\prod_iE_i$ consisting of all $(x_i)_{i\in\N}$ satisfying
$\pi_{ij}(x_j)=x_i$ for all $i\leq j$. Now $\vec\pi(y):=(\pi_i(y))_{i\in\N}$
defines a bijection $\vec\pi:E\to E'$, so there exists a unique measure $\mu$
on the \si-algebra generated by the $(\pi_i(x))_{i\in\N}$ whose image under
$\vec\pi$ equals $\mu'$. By \cite[Lemma~II.18]{Sch73}, this \si-algebra
coincides with the Borel-\si-algebra on $E$.\med

\noi
{\bf Proof of Theorem~\ref{T:max}} In the case without annihilation,
parts~(a)--(e) were proved in \cite[Thm~2]{AS05}. The proof there made
essential use of monotonicity, which is not available in case $a>0$.
Instead of trying to adapt these arguments, replacing monotone convergence by
some other form of convergence wherever necessary, we will make use of
Corollary~\ref{C:allthin}, which will simplify our life considerably.

In view of this, set $\al:=a/(a+c)$ and let $\ov X^{(\infty)}$ be the
$(q,0,(1+\al)b,a+c,\al b+d)$-branco-process started at infinity, as defined in 
\cite[Thm~2]{AS05}. Fix $\eps>0$ and let $(X^\eps_t)_{t\geq\eps}$ be a
$(q,a,b,c,d)$-branco-process started at time $\eps$ in
$X^\eps_\eps=\Thin_{\frac{1}{1+\al}}(\ov X^{(\infty)}_\eps)$. It has been
proved in \cite[Thm~2]{AS05} that $\ov X^{(\infty)}_t\in\Ei_\ga(\La)$ for all
$t\geq 0$ a.s., hence $\Thin_{\frac{1}{1+\al}}(\ov
X^{(\infty)}_\eps)\in\Ei_\ga(\La)$ and $(X^\eps_t)_{t\geq\eps}$ is
well-defined by Proposition~\ref{P:constr}. By Corollary~\ref{C:allthin},
\be
\P[X^\eps_t\in\cdot\,]
=\P\big[\Thin_{\frac{1}{1+\al}}(\ov X^{(\infty)}_t)\in\cdot\,\big]
\qquad(t\geq\eps).
\ee
In particular, this implies that if we construct two processes
$X^\eps,X^{\eps'}$ for two values $0<\eps<\eps'$, then these are consistent in
the sense that $(X^\eps_t)_{t\geq\eps'}$ is equally distributed with
$(X^{\eps'}_t)_{t\geq\eps'}$. By applying the projective limit theorem
sketched above, using the spaces of componentwise cadlag functions from
$(\eps,\infty)$ to $\N^\La$, we may construct a process
$(X^{(\infty)}_t)_{t>0}$ such that $X^{(\infty)}_\eps$ is equally distributed
with $\Thin_{\frac{1}{1+\al}}(\ov X^{(\infty)}_\eps)$ for all $\eps>0$ and
$(X^{(\infty)}_t)_{t>0}$ evolves as a $(q,a,b,c,d)$-branco-process. Let
$\Xc^\dgg$ denote the $(q,r,s,m)$-resem process with $r,s,m$ as in
(\ref{rsmabcd}). Then $\Xc^\dgg$ is dual to both the
$(q,a,b,c,d)$-branco-process (with parameter $\al=a/(a+c)$ in the duality
function) and to the $(q,0,(1+\al)b,a+c,\al b+d)$-branco-process (with duality
function $\Psi(x,\phi)=(1-\phi)^x$).  We have
\be\ba{l}\label{maxdu}
\dis\E\big[(1-(1+\al)\phi)^{\txt X^{(\infty)}_t}\big]
=\E\big[(1-(1+\al)\phi)^{\txt\Thin_{\frac{1}{1+\al}}(\ov X^{(\infty)}_t)}\big]
=\E\big[(1-\phi)^{\txt\ov X^{(\infty)}_t}\big]\\[5pt]
\dis\quad=\P^\phi\big[\Xc^\dgg_t=0]
\qquad(t\geq 0,\ \phi\in[0,1]^\La,\ |\phi|<\infty),
\ec
where the last equality follows from \cite[formula~(5.5)]{AS05} and we assume
$|\phi|<\infty$ to make sure the infinite products are well-defined.
It has been shown in \cite[Thm~2~(d)]{AS05} that the law of $\ov X^{(\infty)}_t$
converges as $t\to\infty$ to an invariant law of the $(q,0,(1+\al)b,a+c,\al
b+d)$-branco-process. By thinning, it follows that the law of $X^{(\infty)}_t$
converges as $t\to\infty$ to an invariant law $\ov\nu$ of the
$(q,a,b,c,d)$-branco-process. Taking the limit $t\to\infty$ in (\ref{maxdu})
we arrive at (\ref{extinct}). Setting
\be
r:=(1+\al)b+(a+c)-(\al b+d)=a+b+c-d,
\ee
we obtain from \cite[Thm~2~(b)]{AS05} and the fact that
$X^{(\infty)}$ is a $1/(1+\al)$-thinning of $\ov X^{(\infty)}$, that
\be
\E[X^{(\infty)}_t(i)]\leq\left\{\ba{cl}
\frac{1}{1+\al}\frac{r}{(a+c)(1-e^{-rt})}\quad&\mbox{if }r\neq 0,\\[5pt]
\frac{1}{1+\al}\frac{1}{(a+c)t}\quad&\mbox{if }r=0\ea\right.
\qquad(i\in\La),
\ee
which by the fact that $1/(1+\al)=(a+c)/(2a+c)$ yields (\ref{explicit}).
Formula~(\ref{thinmax}) is a simple consequence of the way we have defined
$X^{(\infty)}$ as a thinning of $\ov X^{(\infty)}$. This completes the proof
of parts~(a), (b), and (d)--(f) of the theorem.

To prove also part~(c), by formula (\ref{maxdu}) and duality, it suffices to
show that for each $t>0$
\be\label{maxapp}
\E\big[\big(1-(1+\al)\phi\big)^{\txt X^{(n)}_t}\big]
=\E^\phi\big[\big(1-(1+\al)\Xc^\dgg_t\big)^{\txt x^{(n)}}\big]
\asto{n}\P\big[\Xc^\dgg_t=0]\qquad(\phi\in[0,1]^\La,\ |\phi|<\infty).
\ee
By Lemma~\ref{L:notone}~(i) below, $\Xc^\dgg_t(i)<1$ a.s.\ for all $i\in\La$,
hence a.s.\ on the event $\Xc^\dgg_t\neq 0$ there exists some $i\in\La$ such
that $0<\Xc^\dgg_t(i)<1$. It follows that $|1-(1+\al)\Xc^\dgg_t|^{x^{(n)}}\to
0$ as $n\to\infty$ a.s.\ on the event that $\Xc^\dgg_t\neq 0$, hence
(\ref{maxapp}) follows from bounded pointwise convergence.\qed

\noi
{\bf Remark} Let $X^{(n)}$ be as in Theorem~\ref{T:max}~(c). Then, using
Proposition~\ref{P:subdup}, copying the proof of \cite[Thm~2~(b)]{AS05}, we
obtain the uniform estimate
\be\label{explicit2}
\E[X^{(\infty)}_t(i)]\leq
\left\{\ba{cl}
\dis\frac{r'}{(2a+c)(1-e^{-r't})}\quad&\mbox{if }r'\neq 0,\\[8pt]
\dis\frac{1}{(2a+c)t}\quad&\mbox{if }r'=0\ea\right.
\qquad(i\in\La),
\ee
where $r':=2a+b+c-d$. It is easy to see that this estimate is always worse
than the estimate (\ref{explicit}) that we obtained with the help of thinning
(Corollary~\ref{C:allthin}).\med

\subsection{Particles everywhere}

The aim of this section is to prove Lemma~\ref{L:ewhere} below, which, roughly
speaking, says that if we start a $(q,a,b,c,d)$-branco-process in a nontrivial
spatially homogeneous initial law, then for each $t>0$, if we look at
sufficiently many sites, then we are sure to find a particle somewhere. For
zero annihilation rate, this has been proved in \cite[Lemma~6]{AS05}. Results
of this type are well-known, see e.g. the proof of
\cite[Thm~III.5.18]{Lig85}. It seems the main idea of the proof, and in
particular the use of H\"older's inequality in (\ref{plusR}) below or in
\cite[(III.5.30)]{Lig85} goes back to Harris \cite{Har76}. Another essential
ingredient of the proof is some form of almost independence for events that
are sufficiently far apart. For systems where the number of particle per site
is bounded from above, such asymptotic independence follows from
\cite[Thm~I.4.6]{Lig85}, but for branco-processes, the uniform estimate given
there is not available. In \cite{AS05}, we solved this problem by using
monotonicity, which is also not available in the presence of annihilation. 
Instead, we will base our proof on the covariance estimate from
Proposition~\ref{P:covest} above.

\bl{\bf(Particles at the origin)}\label{L:at0}
Let $G$ be a transitive subgroup of ${\rm Aut}(\La)$ and let $\mu$ be a
$G$-homogeneous probability measure on $\Ei_\ga(\La)$. Assume that
$b>0$. Then, for a.e.\ $x$ w.r.t.\ $\mu$, the $(q,a,b,c,d)$-branco-process
started in $X_0=x$ satisfies
\be
\P^x[X_t(0)>0] >0 \qquad(t>0).
\ee
\el
{\bf Proof} Although the statement is intuitively obvious, some work is needed
to make this rigorous. If $a=0$, then by monotonicity (see Lemma~\ref{L:comp},
which extends to infinite initial states by Corollary~\ref{C:finapp}), it
suffices to prove that for a.e.\ $x$ w.r.t.\ $\mu$, there exists some
$i\in\La$ with $x(i)>0$ such that there is a positive probability that a
random walk with jump rates $q$, started in $i$, is at time $t$ in the
origin. Since we are only assuming a weak form of irreducibility (see
(\ref{qassum})~(ii)), this is not entirely obvious, but it is nevertheless
true as has been proved in \cite[Lemma~31]{AS05}.

If $a>0$, then, to avoid problems stemming from the non-monotonicity of
$X$, we use duality. Let $\al,r,s,m$ be as in (\ref{rsmabcd}) and
observe that $m>0$ by our assumptions that $a,b>0$. Define
$\de_0\in[0,1]^{\La}$ by $\de_0(i):=1_{\{i=0\}}$. Then, by duality
(Proposition~\ref{P:dual}), letting $\Xc$ denote the
$(q^\dgg,r,s,m)$-resem-process started in $\Xc_0=\de_0$, we have
\be
\E^x\big[(1-(1+\al))^{\txt X_t(0)}\big]
=\E^{\de_0}\big[(1-(1+\al)\Xc^\dgg_t)^x\big],
\ee
and our claim will follow once we show that for all $t>0$, this quantity is
strictly less than one for a.e.\ $x$ w.r.t.\ $\mu$. Thus, it suffices to show
that $\P^{\de_0}[0<\Xc_t(i)<1]>0$ for some $i\in\La$ such that $x(i)>0$. By
the fact that $m>0$ and Lemma~\ref{L:notone}~(i) below, this can be relaxed to
showing that $\P^{\de_0}[\Xc_t(i)>0]>0$ for some $i\in\La$ such that $x(i)>0$.
Letting $\ti X$ denote the $(q,0,s,r,m)$-branco-process, using duality again
(this time with $\al=0$), it suffices to show that
\be
1>\E^{\de_0}\big[(1-\Xc^\dgg_t)^x\big]=
\E^x\big[0^{\txt\,\ti X_t(0)}\big]=\P^x\big[\ti X_t(0)=0\big].
\ee
Thus, the statement for systems with annihilation rate $a>0$ follows from the
statement for systems with $a=0$.\qed

\bl{\bf(Finiteness of moments)}\label{L:finsec}
Let $X$ be a $(q,a,b,c,d)$-branco-process started in an arbitrary initial law
on $\Ei_\ga(\La)$. Assume that $(\La,q)$ is homogeneous and that $a+c>0$. Then
\be
\E[X_t(i)^m]<\infty\qquad(m\geq 1,\ i\in\La,\ t>0).
\ee
\el
{\bf Proof} By Lemma~\ref{L:comp} and Corollary~\ref{C:finapp}, for each $t>0$
we can couple a $(q,a,b,c,d)$-branco-process $X$ started in an arbitrary
initial law on $\Ei_\ga(\La)$ to the $(q,0,b,a+c,d)$-branco-process $X'$
started in the same initial law, in such a way that $X_t\leq X'_t$ a.s. In
view of this, it suffices to prove the statement for the system $X'$ with zero
annihilation rate and annihilation rate $c':=a+c$. Let $X^{'(n)}$ be the
$(q,0,b,c',d)$-branco-process started in $X^{'(n)}(i)=X'_0(i)\vee n$
$(i\in\La)$. Then, by \cite[Theorem~2~(c)]{AS05}, for each $t>0$ the process
$X^{'(n)}_t$ can be coupled to the process started at infinity, denoted by
$X^{(\infty)}$, in such a way that $X^{'(n)}_t\up X^{(\infty)}_t$ a.s. In view
of this, it suffices to prove that for the process without annihiation started
at infinity
\be
\E[X^{(\infty)}_t(i)^m]<\infty\qquad(m\geq 1,\ i\in\La,\ t>0).
\ee
Let $X^{(n)}$ denote the $(q,0,b,c',d)$-branco-process started in the constant
initial state $X^{(n)}_0(i)=n$ $(i\in\La)$. Again by
\cite[Theorem~2~(c)]{AS05}, it suffices to find upper bounds on
$\E[X^{(n)}_t(i)^m]$ that are uniform in $n$. Such upper bounds have been
derived in \cite[Lemma~(2.13)]{DDL90} for branching-coalescing particle
systems on $\Z^d$ with more general branching mechanisms than considered in
the present paper. In particular, their result includes
$(q,0,b,c',d)$-branco-processes on $\Z^d$ with $c'>0$. Their arguments are
not restricted to $\Z^d$ and apply more generally to underlying lattices
$\La$ and jump kernels $q$ as considered in the present paper, as long as
$(\La,q)$ is homogeneous.\qed

\noi
{\bf Remark} It seems likely that the assumption in Lemma~\ref{L:finsec} that
$(\La,q)$ is homogeneous is not needed. The proof of
\cite[Lemma~(2.13)]{DDL90}, which we apply here, uses translation invariance
in an essential way, however. Since we do not need Lemma~\ref{L:finsec} in the
inhomogeneous case, we will be satisfied with the present statement. It does
not seem easy to adapt the proof of formula (\ref{explicit}) (which holds
without a homogeneity assumption) to obtain estimates for higher moments.

\bl{\bf(Systems with particles everywhere)}\label{L:ewhere}
Assume that $(\La,q)$ is infinite and homogeneous, $G$ is a transitive
subgroup of ${\rm Aut}(\La,q)$, and $a+c>0$, $b>0$. Let $X$ be a
$(q,a,b,c,d)$-branco-process started in a $G$-homogeneous nontrivial initial
law on $\Ei_\ga(\La)$. Then, for any $t>0$,
\be\label{ewhere}
\lim_{n\to\infty}\P[\Thin_{\phi_n}(X_t)=0]=0
\ee
for all $\phi_n\in[0,1]^\La$ satisfying $|\phi_n|\to\infty$.
\el
{\bf Proof} By Lemma~\ref{L:finsec}, restarting the process at some small
positive time if necessary, we can without loss of generality assume that
$\E[X_0(0)^2]<\infty$. Set $\pi_n:=\phi_n/|\phi_n|$ and let $\P^x$ denote the
law of the process started in a deterministic initial state $x$. Then, for
each $r<\infty$ and $t>0$, we can choose $n$ sufficiently large such that
$r\leq|\phi_n|$. Then a $r\pi_n$-thinning is stochastically less than a
$\phi_n$-thinning and therefore
\be\ba{l}\label{prev}
\dis\P^x[\Thin_{\phi_n}(X_t)=0]
\leq\P^x[\Thin_{r\pi_n}(X_t)=0]\\[5pt]
\dis\quad=\E^x\big[\prod_i
(1-r\pi_n(i))^{\txt X_t(i)}\big]
\leq\dis\E^x\big[\prod_i\ex{-r\sum_i\pi_n(i)X_t(i)}\big]\\[5pt]
\dis\quad=:\prod_{i\in A_n}\E^x\big[
\ex{\txt-r\pi_n(i)X_t(i)}\big]+R_n(x)
\leq\prod_{i\in A_n}\E^x\big[\ex{\txt-X_t(i)}\big]^{r\pi_n(i)}+R_n(x),
\ec
where in the last step we have applied Jensen's inequality to the concave
function $z\mapsto z^{r\pi_n(i)}$. For the process started in a nontrivial
homogeneous initial law, we obtain, using H\"older's inequality, for all $n$
sufficiently large such that $r\leq|\phi_n|$,
\bc\label{plusR}
\dis\P[\Thin_{\phi_n}(X_t)=0]
&=&\dis\int\P[X_0\in\di x]\P^x[\Thin_{\phi_n}(X_t)=0]\\[5pt]
&\leq&\dis\int\P[X_0\in\di x]\Big[\prod_{i\in A_n}
\E^x\big[\ex{\txt-X_t(i)}\big]^{r\pi_n(i)}+R_n(x)\Big]\\[5pt]
&\leq&\dis\prod_{i\in A_n}\Big(\int\P[X_0\in\di x]
\E^x\big[\ex{\txt-X_t(i)}\big]^r\Big)^{\pi_n(i)}+\E[R_n(X_0)]\\[5pt]
&=&\dis\prod_{i\in A_n}\Big(\int\P[X_0\in\di x]
\E^x\big[\ex{\txt-X_t(0)}\big]^r\Big)^{\pi_n(i)}+\E[R_n(X_0)]\\[5pt]
&=&\dis\int\P[X_0\in\di x]\E^x\big[\ex{\txt-X_t(0)}\big]^r+\E[R_n(X_0)],
\ec
where we have used spatial homogeneity in the last step but one.

By Corollary~\ref{C:exp}, the quantity $R_n(x)$ defined in (\ref{prev})
can be estimated as
\be
|R_n(x)|\leq\ffrac{1}{2}r^2\sum\subb{k,l}{k\neq l}
\Big(\sum_ix(i)K_t(i;k,l)+\sum_{i,j}x(i)x(j)L_t(i,j;k,l)\Big)\pi_n(k)\pi_n(l).
\ee
It follows that
\bc\label{ER}
\dis\E\big[\big|R_n(X_0(0))\big|\big]
&\leq&\dis\ffrac{1}{2}r^2\sum\subb{k,l}{k\neq l}
\Big(\sum_i\E[X_0(i)]K_t(i;k,l)
+\sum_{i,j}\E[X_0(i)X_0(j)]L_t(i,j;k,l)\Big)\pi_n(k)\pi_n(l)\\[5pt]
&\leq&\dis\ffrac{1}{2}r^2\sum\subb{k,l}{k\neq l}
\Big(\E[X_0(0)]\sum_iK_t(i;k,l)
+\E[X_0(0)^2]\sum_{i,j}L_t(i,j;k,l)\Big)\pi_n(k)\pi_n(l)\\[5pt]
&=:&\dis r^2\sum_{k,l}C(k,l)\pi_n(k)\pi_n(l),
\ec
where by definition $C(k,k):=0$ and we have used that by Cauchy-Schwartz and
translation invariance:
\be
\big|\E[X_0(i)X_0(j)]\big|\leq\E[X_0(i)^2]^{1/2}\E[X_0(j)^2]^{1/2}=\E[X_0(0)^2].
\ee
We claim that
\be\label{Cpp}
\sum_{k,l}C(k,l)\pi_n(k)\pi_n(l)\asto{n}0.
\ee
To see this, we observe that by (\ref{KLhom}), (\ref{KLsum}) and our
assumption that $\E[X_0(0)^2]<\infty$,
\be\label{Cprop}
C(gk,gl)=C(k,l)\quad(g\in G)
\quand
\sum_kC(k,0)<\infty.
\ee
Since $G$ is transitive, for each $l\in\La$ we can choose some $g_l\in G$ such
that $g_ll=0$. In view of this, (\ref{Cprop}) shows in particular that for
each $\eps>0$, the quantity
\be
\big|\{k\in\La:C(k,l)\geq\eps\}\big|
=\big|\{g_lk\in\La:C(g_lk,0)\geq\eps\}\big|
=\big|\{j\in\La:C(j,0)\geq\eps\}\big|=:K_\eps
\ee
does not depend on $l\in\La$ and is finite. It follows that
\be
\sum_l\pi_n(l)\sum_kC(k,l)\pi_n(k)
\leq\sum_l\pi_n(l)\Big(\sum_{k:\,C(k,l)\geq\eps}\pi_n(k)
+\sum_{k:\,C(k,l)<\eps}\pi_n(k)\Big)
\leq K_\eps/|\phi_n|+\eps.
\ee
Since $|\phi_n|\to\infty$ and $\eps>0$ is arbitrary, this proves
(\ref{Cpp}). By (\ref{plusR}) and (\ref{ER}), we conclude that for each
$r<\infty$,
\be
\limsup_{n\to\infty}\P[\Thin_{\phi_n}(X_t)=0]
\leq\int\P[X_0\in\di x]\E^x\big[\ex{\txt-X_t(0)}\big]^r.
\ee
Letting $r\to\infty$, using $b>0$ and Lemma~\ref{L:at0}, we arrive at
(\ref{ewhere}).\qed

\subsection{Long-time limit law}\label{S:hom}

In this section, we prove Theorem~\ref{T:hom}. We first need some preparatory
results.

\bl{\bf(Not exactly one)}\label{L:notone}
Let $\Xc$ be a $(q,r,s,m)$-resem process started in a finite initial state
$\phi\in[0,1]^\La$, $|\phi|<\infty$. Assume that $(\La,q)$ is infinite and
homogeneous and that $m>0$. Then
\begin{itemize}
\item[{\rm(i)}] $\P^\phi[\Xc_t(i)=1]=0$ for each
$t>0$, $i\in\La$.
\item[{\rm(ii)}] $\P^\phi[0<|\Xc_t\wedge(1-\Xc_t)|<K]\to 0$ as $t\to\infty$
  for all $K<\infty$.
\end{itemize}
\el
{\bf Proof} Let $\Xc^+$ and $\Xc^-$ satisfy $\Xc^+_0=\Xc^-_0=\Xc_0=\phi$ and
be given, for times $t>0$, by the solutions to the stochastic differential
equations
\bc\label{indepdif}
\dis\di\Xc^+_t(i)
&=&\dis\big(|q|+s\big)(1-\Xc^+_t(i))\di t-m\Xc^+_t(i)\di t+\sqrt{2r\Xc^+_t(i)(1-\Xc^+_t(i))}\di B_t(i)\\[5pt]
\dis\di\Xc^-_t(i)
&=&\dis-\big(|q|+m\big)\Xc^-_t(i)\di t
+\sqrt{2r\Xc^-_t(i)(1-\Xc^-_t(i))}\di B_t(i)\qquad(t\geq 0,\ i\in\La),
\ec
where $|q|:=\sum_jq(i,j)$, which does not depend on $j\in\La$ by the
transitivity of ${\rm Aut}(\La,q)$, and $(B(i))_{i\in\La}$ is the same
collection of independent Brownian motions as those driving $\Xc$.
By the arguments used in the proof of \cite[Lemma~18]{AS05}, solutions of
(\ref{indepdif}) are pathwise unique and satisfy
\be\label{XXX}
\Xc^-_t\leq\Xc_t\leq\Xc^+_t\quad(t\geq 0)\quad{\rm a.s.}
\ee
Moreover, since (\ref{indepdif}) contains no interaction terms, the
$[0,1]^2$-valued processes $(\Xc^-_t(i),\Xc^+_t(i))_{t\geq 0}$ are independent
for different values of $i\in\La$. Since $\Xc^+(i)$ is a one-dimensional
diffusion with (by grace of the fact that $m>0$) the drift on the boundary
point 1 pointing inwards, it can be proved by standard methods that
\be\label{notone}
\P[\Xc^+_t(i)=1]=0\qquad(t>0,\ i\in\La).
\ee
We defer a precise proof of this fact to Lemma~\ref{L:A1} in the appendix.
Together with (\ref{XXX}), formula (\ref{notone}) proves part~(i) of the lemma.

To prove also part~(ii), we observe that
\be\label{posexp}
\E[\Xc^-_t(i)]=\ex{-(|q|+m)t}\phi(i)\qquad(t>0,\ i\in\La).
\ee
With a bit of work, it is possible to show that there exists a $t_0>0$
and function $(0,t_0]\ni t\mapsto c_t>0$ such that
\be\label{Ebd}
\E\big[\Xc^-_t(i)\wedge(1-\Xc^+_t(i))\big]\geq c_t\phi(i)
\qquad(0<t\leq t_0).
\ee
A precise proof of this fact can be found in Lemma~\ref{L:A2} of the appendix.
We note that for any $[0,1]$-valued random variable $Z$, one has
$\var(Z)=\E[(Z-\E[Z])^2]\leq\E[|Z-\E[Z]|]\leq\E[Z+\E[Z]]=2\E[Z]$.
Applying this to $Z=\Xc^-_t(i)\wedge(1-\Xc^+_t(i))$, using (\ref{posexp}),
we see that
\be\label{Vbd}
\var\big(\Xc^-_t(i)\wedge(1-\Xc^+_t(i))\big)\leq 2\ex{-(|q|+m)t}\phi(i)
\qquad(t>0,\ i\in\La).
\ee
Now (\ref{Ebd}) implies $\E[|\Xc^-_t\wedge(1-\Xc^+_t)|]\geq c_t|\phi|$, while
by (\ref{Vbd}) and the independence of coordinates $i\in\La$,
\be
\var\big(|\Xc^-_t\wedge(1-\Xc^+_t)|\big)
\leq2\ex{-(|q|+m)t}|\phi|\qquad(0<t\leq t_0).
\ee
Since $\Xc_t(i)\wedge(1-\Xc_t(i))\geq\Xc^-_t(i)\wedge(1-\Xc^+_t(i))$, by
Chebyshev, it follows that
\be\label{mid}
\P^\phi\big[|\Xc_t\wedge(1-\Xc_t)|\leq\ffrac{1}{2}c_t|\phi|\big]
\leq\frac{2e^{-(|q|+m)t}|\phi|}{\frac{1}{4}c_t^2|\phi|^2}
\qquad(0<t\leq t_0),
\ee
which tends to zero for $|\phi|\to\infty$. By \cite[Lemma~5]{AS05},
\be
\P^\phi\big[0<|\Xc_t|<K\big]\asto{t}0\qquad(K<\infty).
\ee
It follows that we can choose $L_t\to\infty$ slow enough such that
\be
\P^\phi\big[0<|\Xc_t|<L_t\big]\asto{t}0.
\ee
By (\ref{mid}), we conclude that
\be\ba{l}
\dis\limsup_{t\to\infty}\P^\phi\big[0<|\Xc_t\wedge(1-\Xc_t)|<K\big]\\[5pt]
\dis\quad\leq\limsup_{t\to\infty}
\P^\phi\big[0<|\Xc_t\wedge(1-\Xc_t)|<K\,\big|\,0<|\Xc_{t-t_0}|<L_{t-t_0}\big]
\P^\phi\big[0<|\Xc_{t-t_0}|<L_{t-t_0}\big]\\[5pt]
\dis\quad\phantom{=}+\limsup_{t\to\infty}
\P^\phi\big[0<|\Xc_t\wedge(1-\Xc_t)|<K\,\big|\,|\Xc_{t-t_0}|\geq L_{t-t_0}\big]
\P^\phi\big[|\Xc_{t-t_0}|\geq L_{t-t_0}\big]\\[5pt]
\dis\quad\leq\limsup_{t\to\infty}
\P^\phi\big[|\Xc_t\wedge(1-\Xc_t)|\leq\ffrac{1}{2}c_{t_0}L_{t-t_0}
\,\big|\,|\Xc_{t-t_0}|\geq L_{t-t_0}\big]\\[5pt]
\dis\quad\leq\limsup_{t\to\infty}
\frac{2e^{-(|q|+m)t}L_{t-t_0}}{\frac{1}{4}c_{t_0}^2L_{t-t_0}^2}=0.
\ec
\qed

\noi
{\bf Remark} It seems likely that the condition $m>0$ in Lemma~\ref{L:notone}
is not necessary, at least for part~(i). Indeed, it seems likely that
$(q,r,s,m)$-resem-processes have the `noncompact support property'
\be\label{infsupport}
\P\big[\Xc_t(i)>0,\ \Xc_t(j)=0\big]=0\qquad\big(t>0,\ i,j\in\La,\ q(i,j)>0\big),
\ee
similar to what is known for super random walks \cite{EP91}. Since proving
(\ref{infsupport}) is quite involved and we don't know a reference, we will be
satisfied with proving Lemma~\ref{L:notone} only for $m>0$, which is
sufficient for our purposes.\med

\bl{\bf(Systems with particles everywhere)}\label{L:omnipres}
Assume that $(\La,q)$ is infinite and homogeneous and that $G$ is a
transitive subgroup of ${\rm Aut}(\La,q)$ and $a+c>0$, $b>0$. Let $X$ be the
$(q,a,b,c,d)$-branco process started in a $G$-homogeneous nontrivial
initial law $\Li(X_0)$. Then, for any $t>0$ and
$0\leq\al\leq 1$ and for any $\eps>0$, there exists a $K<\infty$ such that
\be\label{posthin}
|\phi|<\infty\mbox{ and }|\phi\wedge(1-\phi)|\geq K\mbox{ implies }
\E\big[|1-(1+\al)\phi|^{X_t}\big]\leq\eps.
\ee
\el
{\bf Proof} We start by proving that if $\phi_n\in[0,1]^\La$ satisfy
$|\phi_n|<\infty$ and $|\phi_n\wedge(1-\phi_n)|\to\infty$, then
\be\label{posthin2}
\lim_{n\to\infty}\E\big[|1-(1+\al)\phi_n|^{X_t}\big]=0.
\ee
Set $\psi_n:=\phi_n\wedge(1-\phi_n)$. Then, for each $i\in\La$, we have and
$\psi_n(i)\leq 1-\phi_n(i)\leq 2-(1+\al)\phi_n(i)$ and
$\psi_n(i)\leq\phi_n(i)\leq(1+\al)\phi_n(i)$, from which we see that
\be
\psi_n(i)-1\leq 1-(1+\al)\phi_n(i)\leq 1-\psi_n(i),
\ee
or, in other words, $|1-(1+\al)\phi_n(i)|\leq 1-\psi_n$. It follows that
\be
\big|\E\big[1-(1+\al)\phi_n^{X_t}\big]\big|
\leq\E\big[|1-(1+\al)\phi_n|^{X_t}\big]
\leq\E\big[(1-\psi_n)^{X_t}\big]=\P\big[\Thin_{\psi_n}(X_t)=0\big],
\ee
which tends to zero by Lemma~\ref{L:ewhere} and our assumption that
$|\psi_n|\to\infty$.

Now imagine that the lemma does not hold. Then there exists some $\eps>0$ such
that for all $n\geq 1$ we can choose $\phi_n$ with $|\phi_n|<\infty$ and
$|\phi_n\wedge(1-\phi_n)|\geq n$ such that
$\E\big[|1-(1+\al)\phi_n|^{X_t}\big]>\eps$. Since this contradicts
(\ref{posthin2}), we conclude that the lemma must hold.\qed

\noi
{\bf Proof of Theorem~\ref{T:hom}} For $a=0$ the statement has been proved in
\cite[Thm~4~(a)]{AS05}, so without loss of generality we may assume that $a>0$.
By Theorem~\ref{T:max}~(e), it suffices to show that
\be\label{TBconv}
\E\big[\big(1-(1+\al)\phi\big)^{\txt X_t}\big]
\asto{t}\P^\phi\big[\exists t\geq 0\mbox{ such that }\Xc^\dgg_t=0\big]
\qquad(\phi\in[0,1]^\La,\ |\phi|<\infty),
\ee
where $\al:=a/(a+c)$ and $\Xc^\dgg$ denotes the $(q^\dgg,a+c,(1+\al)b,\al
b+d))$-resem-process started in $\phi$. By duality (Proposition \ref{P:dual}),
for each $t\geq 1$,
\be\label{XaX}
\E\big[\big(1-(1+\al)\phi\big)^{\txt X_t}\big]
=\E\big[\big(1-(1+\al)\Xc^\dgg_{t-1}\big)^{\txt X_1}\big],
\ee
where $\Xc^\dgg$ is independent of $X$ and started in $\Xc^\dgg_0=\phi$.
For each $K<\infty$, we may write
\be\ba{l}\label{threesplit}
\dis\E\big[\big(1-(1+\al)\Xc^\dgg_{t-1}\big)^{\txt X_1}\big]
=\P[\Xc^\dgg_{t-1}=0]\\[5pt]
\dis\quad+\E\big[\big(1-(1+\al)\Xc^\dgg_{t-1}\big)^{\txt X_1}
\,\big|\,0<|\Xc^\dgg_{t-1}\wedge(1-\Xc^\dgg_{t-1})|<K\big]
\P\big[0<|\Xc^\dgg_{t-1}\wedge(1-\Xc^\dgg_{t-1})|<K\big]\\[5pt]
\dis\quad+\E\big[\big(1-(1+\al)\Xc^\dgg_{t-1}\big)^{\txt X_1}
\,\big|\,K\leq|\Xc^\dgg_{t-1}\wedge(1-\Xc^\dgg_{t-1})|\big]
\P\big[K\leq|\Xc^\dgg_{t-1}\wedge(1-\Xc^\dgg_{t-1})|\big].
\ec
Here the first term converges, as $t\to\infty$, to $\P^\phi\big[\exists t\geq
  0\mbox{ such that }\Xc^\dgg_t=0\big]$. Note that $\al>0$ by our assumption
that $a>0$. Assume for the moment that also $b>0$. Then
Lemma~\ref{L:notone}~(ii) tells us that the second term on the right-hand side
of (\ref{threesplit}) tends to zero. By Lemma~\ref{L:omnipres}, for each
$\eps>0$ we can choose $K$ large enough such that the third term is bounded in
absolute value by $\eps$. Putting these things together, we arrive at
(\ref{TBconv}).

If $b=0$, then Lemma~\ref{L:notone}~(ii) is not available, but in this case
$|\Xc^\dgg_t|$ is a supermartingale, hence \cite[Lemma~5]{AS05} tells us that
$\P^\phi\big[\exists t\geq 0\mbox{ such that }\Xc^\dgg_t=0\big]=1$, and the
proof proceeds as above.\qed

\appendix

\section{Some facts about coupled Wright-Fisher diffusions}

The aim of this appendix is to prove two simple facts about (coupled)
Wright-Fisher diffusions. In particular, applying Lemmas~\ref{L:A1} and
\ref{L:A2} to $X=\Xc^+(i)$, $Y=\Xc^-(i)$, $a=|q|+s$, $b=m$ and $c=|q|+m$
yields formulas (\ref{notone}) and (\ref{Ebd}), respectively.

For $a,b\geq 0$ and $r>0$, let $X$ denote the pathwise unique (by \cite{YW71})
$[0,1]$-valued solution to the stochastic differential equation
\be\label{WF}
\di X_t=a(1-X_t)\di t-bX_t\di t+\sqrt{2rX(1-X)}\di B_t,
\ee
where $B$ is standard Brownian motion.

\bl{\bf(No mass on boundary)}\label{L:A1}
If $b>0$, then
\be
\P[X_t=1]=0\qquad(t>0),
\ee
regardless of the initial law.
\el
{\bf Proof} If $a,b>0$, then it is well known that $X$ has a transition
density (see Propositions 3 and 4 in \cite{Pal11} along with the discussion on
page 1183 or \cite{Gri79b, Gri79a}). Consequently $\P[X_t=1]=0$ and hence the
result follows. If $a=0$ but $b>0$, then by standard comparision results (see
\cite[Thm.~6.2]{Bas98} or \cite[Lemma~18]{AS05}), if $Z_0=X_0$ and $Z$ solves
the SDE (\ref{WF}) with $a=b/2$ and $b$ replaced by $b/2$, relative to the
same Brownian motion, then $X_t\leq Z_t$ and hence $\P[X_t=1]\leq\P[Z_t=1]=0$
for all $t>0$.\qed

\bl{\bf(Moment dual)}\label{L:A3}
Let $K=(K_t)_{t\geq 0}$ be a Markov process with state space
$\N\cup\{\infty\}$, where $\infty$ is a trap, and $K$ jumps from states
$k\in\N$ as
\be\ba{ll}
k\mapsto k-1\qquad&\mbox{with rate }ak+rk(k-1),\\[5pt]
k\mapsto\infty\qquad&\mbox{with rate }bk.
\ec
Then
\be
\E^x[X^k_t]=\E^k[x^{K_t}]\qquad(t\geq 0,\ x\in[0,1],\ k\in\N),
\ee
where $x^0:=1$ and $x^\infty:=0$ for all $x\in[0,1]$.
\el
{\bf Proof} Let
\bc
\dis\Gi f(x)&:=&\dis\big[a(1-x)\dif{x}-bx\dif{x}
+rx(1-x)\diff{x}\big]f(x),\\[5pt]
\dis Gf(k)&:=&\dis\big[ak+rk(k-1)\big]\big\{f(k-1)-f(k)\}
+bk\big\{f(\infty)-f(k)\big\}
\ec
be the generators of the processes $X$ and $K$,
respectively, and let $\psi(x,k):=x^k$ be the duality function. Then
\be\label{GiG}
\Gi\psi(\,\cdot\,,k)(x)
=ak(x^{k-1}-x^k)-bkx^k+rk(k-1)(x^{k-1}-x^k)=G\psi(x,,\,\cdot\,)(k),
\ee
where the term with $k(k-1)$ is zero for $k=1$ and both sides of the equation
are zero for $k=0$. The claim now follows from \cite[Thm~7]{AS05} and
\cite{AS09b} and the fact that the expression in (\ref{GiG}) is bounded
uniformly in $x$ and $k$, which guarantees the required integrability.

Although this is not needed for the proof, this duality may be understood as
follows. We can view $X_t$ as the frequency of type-one organisms in a large
population where pairs of organisms are resampled with rate $2r$ and organisms
mutate to type 1 and 0, respectively, with rates $a$ and $b$. Then
$\E[X^k_t]$ is the probability that $k$ organisms, sampled from the population
at time $t$, are all of type one. We can view $K_t$ as the ancestors of these
organism at time zero, where we neglect organisms that due to mutation are
sure to be of type one while on the other hand the state $K_t=\infty$
signifies that due to a mutatation event, at least one of these ancestors is
of type zero.\qed

Now let $X$ be as in (\ref{WF}), let $c\geq 0$, and let $Y$ be given by the
pathwise unique solution to the stochastic differential equation
\be\label{WFY}
\di Y_t=-cY_t\di t+\sqrt{2rY(1-Y)}\di B_t,
\ee
driven by the {\em same} Brownian motion as $X$.

\bl{\bf(Feller property)}\label{L:A4}
Let $(X,Y)$ be given by the pathwise unique solutions of (\ref{WF}) and
(\ref{WFY}), and let
$K_t((x,y),\,\cdot\,):=\P^{(x,y)}[(X_t,Y_t)\in\,\cdot\,]$ denote the
transition probabilities of $(X,Y)$. Then the map $(t,x,y)\mapsto
K_t((x,y),\,\cdot\,)$ from $\half\times[0,1]$ into the probability measures on
$[0,1]^2$ is continuous w.r.t.\ weak convergence of probability measures.
\el
{\bf Proof} It follows from well-known results \cite[Corollary~5.3.4 and
  Theorem~5.3.6]{EK86} that pathwise uniqueness for a stochastic differential
equation implies uniqueness of solutions to the martingale problem for the
associated differential operator, which is in our case given by
\be
A:=a(1-x)\dif{x}-bx\dif{x}+rx(1-x)\diff{x}
-cy\dif{y}+ry(1-y)\diff{y}+2r\sqrt{x(1-x)y(1-y)}\difif{x}{y},
\ee
with domain $\Ci^2[0,1]^2$. Now if $(X^n,Y^n)$ are solutions to this
martingale problem with deterministic initial states $(X^n_0,Y^n_0)=(x_n,y_n)$
converging to some limit $(x,y)\in[0,1]^2$, and $(X,Y)$ denotes the process
started in $(x,y)$, then \cite[Lemma~4.5.1 and Remark~4.5.2]{EK86} imply that 
\be
\P[(X^n_t,Y^n_t)_{t\geq 0}\in\cdot\,]
\Asto{n}\P[(X_t,Y_t)_{t\geq 0}\in\cdot\,],
\ee
where $\Rightarrow$ denotes weak convergence of probability laws on the space
$\Ci_{[0,1]^2}\half$ of continuous functions from $\half$ into $[0,1]^2$,
equipped with the topology of locally uniform convergence. In particular, this
implies the stated continuity of the transition probabilities.\qed

\bl{\bf(Linear estimate)}\label{L:A2}
Assume that $b>0$. Then there exists a $t_0>0$ and function
$(0,t_0]\ni t\mapsto\la_t>0$ such that the process started in
$(X_0,Y_0)=(z,z)$ satisfies
\be\label{lest}
\E^{(z,z)}\big[Y_t\wedge(1-X_t)\big]\geq\la_tz
\qquad(0<t\leq t_0,\ 0\leq z\leq 1).  
\ee
\el
{\bf Proof} We estimate 
\be\ba{l}\label{zz}
\dis\E^{(z,z)}\big[Y_t\wedge(1-X_t)\big]
\geq\E^{(z,z)}\big[Y_t(1-X_t)\big]
\geq\ffrac{1}{2}\E^{(z,z)}\big[Y_t1_{\{X_t\leq\frac{1}{2}\}}\big]\\[5pt]
\dis\quad=\ffrac{1}{2}\big(\E^z[Y_t]-\E^z[1_{\{X_t>\frac{1}{2}\}}]\big)
\geq\ffrac{1}{2}\E^z[Y_t]-2\E^z[X^2_t],
\ec
where the last step we have used that $1_{\{x>\frac{1}{2}\}}\leq 4x^2$. By
Lemma~\ref{L:A3},
\be\label{Yexp}
\E^z[Y_t]=e^{-ct}z\qquad(t\geq 0,\ z\in[0,1]), 
\ee
while by the same lemma
\be\ba{l}
\dis\E^z[X^2_t]=\E^2[z^{K_t}]
\leq\P^2[K_t\leq 1]z+\P^2[K_t=2]z^2\\[5pt]
=(1-e^{-2(a+r)t})z+e^{-2(a+b+r)t}z^2
\leq 2(a+r)tz+z^2.
\ec
Combining this with (\ref{zz}) yields 
\be
\E^{(z,z)}\big[Y_t\wedge(1-X_t)\big]
\geq\big(\ffrac{1}{2}e^{-ct}-4(a+r)t-2z\big)z.  
\ee
Choosing $t_0>0$ and $z_0>0$ small enough, we find that 
\be\label{prest}
\E^{(z,z)}\big[Y_t\wedge(1-X_t)\big]
\geq\ffrac{1}{4}z\qquad(0\leq t\leq t_0,\ 0\leq z\leq z_0).  
\ee 
To extend this to all $z\in[0,1]$, at the cost of assuming that $t>0$ and
replacing the constant $1/4$ by a possibly worse, time-dependent constant
$\la_t$, we observe that by Lemma~\ref{L:A4}, the function $[0,1]\ni
z\mapsto\E^{(z,z)}\big[Y_t\wedge(1-X_t)\big]$ is continuous. Since by
Lemma~\ref{L:A1} and (\ref{Yexp}),
\be
\E^{(z,z)}\big[Y_t\wedge(1-X_t)\big]>0\qquad(t>0,\ z\in(0,1]),
\ee
using continuity, we may estimate $\E^{(z,z)}\big[Y_t\wedge(1-X_t)\big]$
uniformly from below on $[z_0,1]$, which together with (\ref{prest}) yields
(\ref{lest}).\qed

\end{document}